\documentclass[11pt]{amsart}
\usepackage[linktocpage]{hyperref}
\usepackage{amssymb, paralist, xspace, graphicx, url, amscd, euscript, mathrsfs,stmaryrd,epic,eepic,color}

\usepackage[all]{xy}
\usepackage{amsthm}
\SelectTips{cm}{}
\numberwithin{equation}{section}

1
\numberwithin{subsection}{section}

\allowdisplaybreaks[1]

\newtheorem*{namedtheorem}{\theoremname}
\newcommand{\theoremname}{testing}

\newtheorem*{maintheorem}{Theorem}

\newtheorem{theorem}{Theorem}
\newtheorem{proposition}[subsection]{Proposition}
\newtheorem{proposition-definition}[subsection]
{Proposition-Definition}
\newtheorem{corollary}[subsection]{Corollary}
\newtheorem{lemma}[subsection]{Lemma}
\theoremstyle{definition}
\newtheorem{definition}[subsection]{Definition}
\newtheorem*{question}{Question}

\newcommand\cA{\mathcal{A}}
\newcommand\cB{\mathcal{B}}
\newcommand\cC{\mathcal{C}}
\newcommand\cD{\mathcal{D}}

\newcommand\cH{\mathcal{H}}

\newcommand\cL{\mathcal{L}}

\newcommand\cT{\mathcal{T}}

\newcommand\CC{\mathbb{C}}

\newcommand\HH{\mathbb{H}}

\newcommand\RR{\mathbb{R}}

\makeatletter
\def\ExtendSymbol#1#2#3#4#5{\ext@arrow 0099{\arrowfill@#1#2#3}{#4}{#5}}
\def\RightExtendSymbol#1#2#3#4#5{\ext@arrow 0359{\arrowfill@#1#2#3}{#4}{#5}}
\def\LeftExtendSymbol#1#2#3#4#5{\ext@arrow 6095{\arrowfill@#1#2#3}{#4}{#5}}
\makeatother

\newcommand\myLongEqual[2][]{\ExtendSymbol{=}{=}{=}{#1}{#2}}
\theoremstyle{plain}
\theoremstyle{definition}

\begin{document}

\title{Teichm$\ddot{\text{U}}$ller Space Is Totally Geodesic In Goldman Space}

\author{Qiongling Li}

\address{Department of Mathematics\\
Rice University\\
6100 Main Street\\
Houston, TX 77005\\
U.S.A.}

\email{ql4@rice.edu}

\begin{abstract}
We construct a new Riemannian metric on Goldman space $\mathcal{B}(S)$, the space of the equivalence classes of convex projective structures on the surface $S$, and then prove the new metric, as well as the metric of Darvishzadeh and Goldman, restricts to be the Weil-Petersson metric on Teichm$\ddot{\text{u}}$ller space, embedded as a submanifold of Goldman space $\mathcal{B}(S)$. Moreover, Teichm$\ddot{\text{u}}$ller space endowed with the Weil-Petersson metric then is totally geodesic in the Riemannian manifold $\mathcal{B}(S)$.
\end{abstract}
\maketitle

%%%%%%%%%%%%%%%%%%%%%%%%%%%%%%%%%%%%%%%%%%%%%%start section
%
%
%
%
\tableofcontents

\begin{section}{Introduction}
An $\RR P^2$-structure on a surface is a system of coordinate charts in $\RR P^2$ with transition maps in $PGL(3,\RR)$. Moreover, for a convex $\RR P^2$-structure on a smooth manifold $M$, we may write $M=\Omega/{\Gamma}$, with $\Omega$ a convex domain in some $\RR^2\subset \RR P^2$ and $\Gamma\subset PGL(3,\RR)$. When $M$ is a closed surface $S$ of genus $g>1$, then the equivalence classes of such structures form a moduli space $\mathcal{B}(S)$ homeomorphic to an open cell of dimension 16$(g-1)$ (see \cite{Goldman 3}). 

Labourie \cite{Labourie} and Loftin \cite{Loftin 2} independently gave the correspondence between the deformation space $\mathcal{B}(S)$ and the space of pairs $(\Sigma, U)$, where $\Sigma$ is a Riemann surface varying in Teichm$\ddot{\text{u}}$ller space and $U$ is a cubic differential on $\Sigma$. Teichm$\ddot{\text{u}}$ller space $\cT(S)$ embeds inside $\mathcal{B}(S)$ as the locus of pairs $(\Sigma, 0)$, where 0 represent the vanishing cubic differentials on the Riemann surface $\Sigma$.

It is of interest to know what of that rich geometric structure extends to $\mathcal{B}(S)$. In \cite{Goldman 4}, a symplectic structure on $\mathcal{B}(S)$ is defined, which extends the symplectic structure on the Teichm$\ddot{\text{u}}$ller space $\cT(S)$ defined by the Weil-Petersson K$\ddot{\text{a}}$hler form. For the Riemannian metric, it is natural to ask the following questions:
\begin{question}(i) Does there exist a Riemannian metric on the deformation space $\mathcal{B}(S)$ restricting to be the Weil-Petersson metric on Teichm$\ddot{\text{u}}$ller space $\cT(S)$?\\
(ii) If (i) is satisfied, is the Teichm$\ddot{\text{u}}$ller space $\cT(S)$ endowed with the Weil-Petersson metric totally geodesic within the deformation space $\mathcal{B}(S)$ endowed with the new metric?
\end{question}

To answer the above questions, we first consider the Riemannian metric on $\mathcal{B}(S)$ constructed by Darvishzadeh and Goldman (see \cite{Goldman 2}), which will be referred to as the DG metric. 

We show that the DG metric answers (i) affirmatively in this paper (Theorem \ref{restriction}). By the nature of Koszul-Vinberg metric, we are not able to see directly whether the DG metric satisfies (ii). 

To address this issue, we make use of the Cheng-Yau metric, which is closely related to the correspondence between $\cB(S)$ and the space of pairs $(\Sigma, U)$, to construct a new Riemannian metric on the deformation space $\mathcal{B}(S)$, will be referred to as the Loftin metric. In fact, Loftin mentioned in \cite{Loftin 3} that the construction of the DG metric can be carried out with other invariant affine K$\ddot{\text{a}}$hler metrics instead of the Koszul-Vinberg metric, e.g., the Cheng-Yau metric.  

Now we can answer both parts (i) and (ii) of the above question affirmatively with the Loftin metric, namely,
\begin{maintheorem}\ref{WP} The DG metric and the Loftin metric both restrict to a constant multiple of the Weil-Petersson metric on Teichm$\ddot{\text{u}}$ller space.
\end{maintheorem}

\begin{maintheorem}\ref{total} Teichm$\ddot{\text{u}}$ller space endowed with the Weil-Petersson metric is totally geodesic in $\mathcal{B}(S)$ endowed with the Loftin metric.
\end{maintheorem}

Recently, M. Bridgeman, D. Canary, F. Labourie and A. Sambarino in \cite{Bridgeman} construct a mapping class group invariant Riemannian metric on the Hitchin component $\cH(S)$ of $Hom(\pi(S), PSL(n,\RR))/PSL(n,\RR)$, which is called pressure metric. They showed that the pressure metric is an extension of the Weil-Petersson metric on the Fuchsian representations from thermodynamical formalism. When restricted to the $SL(3,\RR)$ case, the Hitchin component coincides with the deformation space $\mathcal{B}(S)$ (see \cite{Goldman}), hence the pressure metric also answers part (i) of the question affirmatively.

\begin{subsection}*{Outline of the proof} Firstly, to show Theorem \ref{WP}, we begin with showing that the Loftin metric and the DG metric are isometric (up to a constant multiple) when restricted to the Teichm$\ddot{\text{u}}$ller locus (mainly because the ingredients in the definition of two metrics coincide for the hyperbolic structure case). Then it is sufficient to show that the Loftin metric restricts to be the Weil-Petersson metric on the Teichm$\ddot{\text{u}}$ller locus (Proposition \ref{Loftin metric}). Here, these  two metrics are defined on different descriptions of Teichm$\ddot{\text{u}}$ller space: The Weil-Petersson metric is defined on the usual Teichm$\ddot{\text{u}}$ller space, and the Loftin metric is defined on the Teichm$\ddot{\text{u}}$ller locus in $\cB(S)$. Hence we need to identify tangent vectors of Teichm$\ddot{\text{u}}$ller space and those of the Teichm$\ddot{\text{u}}$ller locus in $\cB(S)$ (Lemma \ref{harmonic} (i)). Moreover, we calculate that the 1-forms taking values in flat $sl(3,\RR)$-bundle we choose are, in fact, harmonic in their cohomology class (harmonicity is essential in the definition of the Loftin metric) (Lemma \ref{harmonic} (ii)). The explicit expression of the metric on the flat $sl(3,\RR)$-bundle in Lemma \ref{pairing g} is extremely helpful for the calculation. Then we finish the proof of Proposition \ref{Loftin metric} by comparing the Loftin pairing of the harmonic representatives and the Weil-Petersson pairing of the original tangent vectors.

Secondly, to show Theorem \ref{total}, we apply a result from Riemannian geometry that the fixed set of an isometry is a totally geodesical submanifold in the original manifold. The remaining goal is to find an automorphism of the deformation space $\mathcal{B}(S)$ satisfying that (1) it has exactly Teichm$\ddot{\text{u}}$ller locus as the fixed set, and (2) it is an isometry with respect to the Loftin metric. We construct a dual map $\tau$ of $\mathcal{B}(S)$ and then show that the dual map $\tau$ satisfies both (1) and (2) to finish the proof of Theorem \ref{total}. Statement (1) immediately follows from the construction of the dual map. In fact, if we use the correspondence of $\mathcal{B}(S)$ and space of pairs $(\Sigma, U)$, this dual map $\tau$ takes $(\Sigma, U)$ to $(\Sigma, -U)$. The fixed set of the dual map $\tau$ is the set $\{(\Sigma,0)\}$, which is exactly the embedding image of Teichm$\ddot{\text{u}}$ller space. The remaining part is to show Statement (2) (Theorem \ref{Z}), namely, the dual map $\tau$ is an isometry. 
\end{subsection}

\begin{subsection}*{Plan of the paper} We organize this paper as follows: In \S \ref{setting}, we introduce convex $\RR P^2$-structures on surfaces and the deformation space $\mathcal{B}(S)$ of convex $\RR P^2$-structures on surfaces.  In \S \ref{correspondence}, we describe the correspondence between the deformation space $\mathcal{B}(S)$ and the space of pairs $(\Sigma, U)$, where $\Sigma$ is a Riemann surface varying in Teichm$\ddot{\text{u}}$ller space and $U$ is a holomorphic cubic differential on $\Sigma$. In \S \ref{new metric}, we construct a new Riemannian metric (called the Loftin metric) on the deformation space $\mathcal{B}(S)$ and give the detail of the construction of the DG metric. In \S \ref{e}, we describe the embedding of Teichm$\ddot{\text{u}}$ller space into $\mathcal{B}(S)$ from algebraic and geometric viewpoints. We devote \S \ref{restriction} to showing that the two metrics defined in \S\ref{new metric} on the deformation space $\mathcal{B}(S)$ both restrict to be the Weil-Petersson metric on Teichm$\ddot{\text{u}}$ller space $\cT(S)$. In \S \ref{geodesic}, we introduce a dual map of the deformation space $\mathcal{B}(S)$ and show that it is an isometry, and then finally prove Theorem \ref{total}.
\end{subsection}

\begin{subsection}*{Acknowledgements} The author wishes to express her gratitude to her advisor Michael Wolf, for suggesting the problem, many helpful discussions and encouragement. Parts of this paper were completed while the author participated events organized by GEAR (Junior Retreat and Retreat 2012) and she wishes to thank the organizers for the hospitality and especially for many fruitful discussions during her visit. The author also wants to thank Song Dai for several interesting conversations.
\end{subsection}
\end{section}

\begin{section}{Deformation space of convex $\RR P^2$-structures}\label{setting}
Let $M$ be a smooth 2-manifold. 
\begin{definition} A real projective structure on $M$ is an atlas of charts $\{(U_\alpha, \psi_{\alpha})\}$, such that\\
(i) $\{U_{\alpha}\}$ is an open cover of $M$;\\
(ii) For each $\alpha$, the map $\psi_{\alpha}:U_{\alpha}\rightarrow \RR P^2$ is a diffeomorphism onto its image; and\\
(iii) The change of coordinates are locally projective: If $\{(U_{\alpha}, \psi_{\alpha})\}$ and $\{(U_{\beta}, \psi_{\beta})\}$ are two such coordinate charts, then the restriction of $\psi_{\beta}\circ{\psi_{\alpha}}^{-1}$ to any connected component of ${\psi_{\beta}}^{-1}(\psi_{\alpha}(U_{\alpha}\cap U_{\beta}))$ is a projective transformation.
\end{definition}
\begin{definition}
A manifold with an $\RR P^2$-structure is called an $\RR P^2$-manifold.
\end{definition}
\begin{definition} An $\RR P^2$-structure on $M$ is called convex if its developing map is a diffeomorphism of $\widetilde{M}$ onto a convex domain $\Omega$  in some affine ${\RR}^2\subset\RR P^2$.  In this case, we can realize $M=\Omega/{\Gamma}$, where $\Gamma$ is a subgroup of $PGL(3,\RR)$ which acts discretely and properly discontinuous on $\Omega$.  Moreover, a convex $\RR P^2$-structure on $M$ is called properly convex if $\Omega$ is bounded.
\end{definition}

Let $S$ be a closed surface of genus $g>1$. Define $\mathcal{B}(S)=\{(f,M)|f:S\rightarrow M$ is a diffeomorphism and $M$ is a convex $\RR P^2$-manifold$\}/{\sim}$, which we refer to as Goldman space. The equivalence relation $\sim$ means, two elements $(f,M),(f^{\prime},M^{\prime})$ are equivalent if and only if there exists a projective isomorphism $h:M\rightarrow M^{\prime}$ such that $h\circ f$ is isotopic to $f^{\prime}$. We have that Goldman space $\mathcal{B}(S)$ is open, and the holonomy map is an embedding of $\mathcal{B}(S)$ to $Hom(\pi,PGL(3,\RR))/PGL(3,\RR)$.

The Zariski tangent space to $Hom(\pi,PGL(3,\RR))/PGL(3,\RR)$ at $[\rho]$ (hence also the tangent space to $\mathcal{B}(S)$ at $[\rho]$) is isomorphic to $H^1(\pi, sl(3,\RR))$ which by de Rham's theorem is isomorphic to $H^1(S; sl(3,\RR)_{Ad\rho})$, where $sl(3,\RR)_{Ad\rho}$ is the flat $sl(3,\RR)$-bundle over $S$ with holonomy representation $Ad\rho$ (see \cite{Goldman 1}, pp. 208-209). Explicitly,  $sl(3,\RR) _{Ad\rho}$ is identified with \begin{equation}\label{bundle}\widetilde{S}\times sl(3,\RR)/\{(\widetilde{s},x)\sim(\gamma \widetilde{s},Ad\rho(\gamma)(x))\}, \end{equation} for all $\gamma\in \pi, \widetilde{s}\in \widetilde{S}, x\in sl(3,\RR)$.
\end{section}

\begin{section}{Correspondence of $\mathcal{B}(S)$ and spaces of pairs $(\Sigma, U)$}\label{correspondence} 
As we mentioned in the introduction, we have another description of Goldman space $\mathcal{B}(S)$ as follows:
\begin{proposition}(Theorem 2 in Loftin \cite{Loftin 2}, Theorem 1.0.2 in Labourie \cite{Labourie}) \label{cubic differential} There exists a natural bijective correspondence between convex $\RR
P^2$-structures on $S$ and pairs $(\Sigma,U)$, where $\Sigma$ is a Riemann surface homeomorphic to $S$, and $U$ is a holomorphic cubic differential on $\Sigma$.
\end{proposition}
Since we rely heavily on the construction of the bijection in Proposition \ref{cubic differential}, we give a version of the arguments here for reader's convenience. The arguments mainly follow Loftin \cite{Loftin 4}. 

Before explaining the detail of the correpondence, we first state the main idea as follows: To start with, given a convex ${\RR P}^2$- structure on the surface $S$, we write $M\cong\Omega/{\Gamma}$, where $\Omega$ is a bounded convex domain in $\RR^2$. For a bounded convex domain $\Omega$, there is a unique hypersurface asymptotic to the boundary of the open cone $\cC\subset {\RR}^3$ above $\Omega$ called the hyperbolic affine sphere (which will be defined later) (see Proposition \ref{affine sphere}). This hyperbolic affine sphere $H\subset \cC$ is invariant under automorphisms of $\cC$ in $SL(3,\RR)$. The restriction of the projection map $\pi:\cC\rightarrow \Omega$ induces a diffeomorphism of $H$ onto $\Omega$. Affine differential geometry provides an $SL(3,\RR)$-invariant structure on the hyperbolic affine sphere $H$ which then descends to $M=\Omega/{\Gamma}$. Then the affine metric on the surface (will be defined later) induces a conformal structure, hence gives a Riemann surface structure $\Sigma$ on the surface. Moreover, the difference of the Levi-Civita connection of the affine metric on $H$ and the Blaschke connection of $H$ (will be defined later) induces a holomorphic cubic differential on the Riemann surface $\Sigma$.

\begin{subsection}*{Hyperbolic Affine Sphere} Consider a hypersurface immersion $f:H\rightarrow {\RR}^3$, and consider a transversal vector field $\xi$ on the hypersurface $H$. We have the equations:
\begin{eqnarray*}
D_X f_*(Y)&=&f_*(\nabla_XY)+h(X,Y)\xi\\
D_X{\xi}&=&-f_*(SX)+\beta(X)\xi.
\end{eqnarray*}
Here, $X$ and $Y$ are tangent vectors on $H$, the operator $D$ is the canonical flat connection induced from ${\RR}^3$, the operator $\nabla$ is a torsion-free connection, the form $h$ is a symmetric bilinear form on $T_x H$, the map $S$ is an endormorphism of $T_x H$, and $\beta$ is a one-form. An affine normal of $H$ is a transversal vector field which is invariant under affine automorphisms of $H$. 

An affine sphere is a hypersurface $H$ in ${\RR}^3$ satisfying the condition that all its affine normals point toward a given point in ${\RR}^3$, called the center. Moreover, if the center lies on the concave side of $H$; and if the map $S=LI$, where the affine mean curvature $L$ is a constant negative function on $H$ and $I$ is the identity map, we call $H$ is a hyperbolic affine sphere.

Thus by scaling, we can normalize any hyperbolic affine sphere to have $L=-1$. Also, we can translate so that the center is 0, i.e., from now on, we restrict to hyperbolic affine spheres with center at origin and affine mean curvature -1. In this case, the affine normal $\xi=f$, where $f$ is the embedding of $H$ into ${\RR}^3$. The structure equations (see \cite{Loftin 4}) then become 
\begin{eqnarray*}
D_X f_*(Y)&=&f_*(\nabla_XY)+g(X,Y)f\\
D_X{f}&=&f_*(X).
\end{eqnarray*}
\end{subsection}The connection $\nabla$ is called the Blaschke connection. The bilinear form $g$ is called the Blaschke metric, or the affine metric.

\begin{proposition}(Cheng and Yau \cite{Cheng-Yau},\cite{Cheng-Yau 1}, Calabi and Nirenberg (with clarification by Gigena \cite{Gigena}, Sasaki \cite{Sataki} and A.M Li \cite{Li},\cite{Li1})) \label{affine sphere}
Consider a convex, bounded domain $\Omega\subset {\RR}^2$, where ${\RR}^2$ is embedded in ${\RR}^3$ as the affine space $\{x_3=1\}$. Then, there is a unique properly embedded hyperbolic affine sphere $H\subset {\RR}^3$ of affine mean curvature -1 and center 0 asymptotic to the boundary of the cone $\cC\subset {\RR}^3$.
\end{proposition}

\begin{subsection}*{The tautological bundle} We define ${\RR P}^2$ as the space of all lines $l$ passing through 0 in ${\RR}^3$. Then the subset of ${\RR P}^2\times {\RR}^3$ consisting of all $(p,l)$ with $p\in l$ is the total space for the tautological line bundle $\cL$ of ${\RR P}^2$. Given an ${\RR P}^2$-manifold $M$, the bundle ${dev}^{-1}\cL$ defines the tautological bundle on $\widetilde{M}$. We say $\widetilde{M}$ admits a tautological bundle if this structure descends to $M$, i.e., if there is a line bundle on $M$ which pulls back to ${dev}^{-1}\cL$ on $\widetilde{M}$ under the universal covering map. For simplicity, we denote this line bundle as $\cL$ also. By Proposition 2.2.1 in Loftin \cite{Loftin 2}, a manifold $M$ with convex ${\RR P}^2$- structure admits an oriented tautological bundle.
\end{subsection}

\begin{subsection}*{Affine sphere structure}Let $M$ be an ${\RR P}^2$-manifold with oriented tautological bundle $\cL$. Then the total space of the positive part of $\cL$ (i.e., the ${\RR}^+$ part of each fiber of the line bundle $\cL$) is locally a cone in ${\RR}^3$. We say $M$ admits an affine sphere structure if there is a section $s$ of $\cL$ so that for each coordinate chart $\mathcal{U}$ of $M$, $s(\mathcal{U})$ is a hyperbolic affine sphere with center 0 and affine mean curvature -1 in the cone $\cC$. \\
\end{subsection}

Combining Proposition 2.2.1 and Theorem 4 in Loftin \cite{Loftin 2} with the fact that any convex $\RR P^2$-structure on a compact surface $S$ must be properly convex by Kuiper \cite{Kuiper}, we have the following proposition:
\begin{proposition}\label{affine}  Let $M$ be a convex $\RR P^2$-manifold homeomorphic to $S$, then we have\\
1. $M$ admits a negative strictly convex section $u$ of the dual tautological bundle $\cL^*$ satisfying $det(u_{ij})=(\frac{1}{u})^{4}$ so that the metric $\frac{-u_{ij}}{u}$ is complete.\\
2. $M$ admits an affine sphere structure whose metric is complete.
\end{proposition}

Now suppose we have $M\cong\Omega/{\Gamma}$; by Proposition \ref{affine sphere}, we have a unique hyperbolic affine sphere $H$ in the cone. Consider a local conformal coordinate $z=x+iy$ on the hyperbolic affine sphere $H$ with center the origin and affine curvature -1. Then the affine metric is given by $g=e^{\psi}|dz|^2$ for some function $\psi$. Parametrize the surface by \begin{equation}\label{parametrization}
f:\mathcal{D}\rightarrow {\RR}^3,~~~~\text{with $\mathcal{D}$ a domain in $\CC$.}
\end{equation} Then we have the following structure equations for the affine sphere:
\begin{eqnarray}
D_XY&=&\nabla_XY+g(X,Y)f\label{metric splitting}\\
D_Xf&=&X\nonumber
\end{eqnarray}
Here $D$ is the canonical flat connection on ${\RR}^3$, $\nabla$ is a projectively flat connection, and $g$ is the affine metric.

We also consider $\widehat{\nabla}$ the Levi-Civita connection with respect to the affine metric $g$. To get a holomorphic cubic differential from this construction, we consider the Pick form $C:=\widehat{\nabla}-\nabla$ which is a tensor measuring the difference between the Levi-Civita connection and the Blaschke connection. In index notation, we have the following conditions (see Theorem 4.3 in \cite{Nomizu}) \begin{equation*}
\Sigma_i C_{ij}^i=0 ~~~\text{for all} j; ~~~~~C_{ijk}~~~~\text{symmetric in}~i,j,k,\end{equation*}where we use $g$ to lower the index. In addition, if $C$ vanishes identically on the hyperbolic affine sphere $H$ with center the origin and affine curvature -1, then $H$ must be the hyperboloid in ${\RR}^3$ (see Theorem 4.5 in \cite{Nomizu}). The symmetries of the Pick form show that it has only two linearly independent factors, which are realized as the real and complex parts of a holomorphic cubic differential $U$ on $S$ under the complex structure with complex $z$ coordinate.
\end{section}
\begin{section}{A new Riemannian metric on $\mathcal{B}(S)$}\label{new metric}
Darvishzadeh and Goldman \cite{Goldman 2} construct a Riemannian metric, which will be referred to as the DG metric, on the deformation space $\mathcal{B}(S)$ using the Koszul-Vinberg metric on the cone. 
In this section, we first give the construction of a new Riemannian metric, will be referred to as the Loftin metric, defined on the deformation space $\mathcal{B}(S)$ but using the Cheng-Yau metric (defined below) on the cone and then give the construction of the DG metric. 

\begin{subsection}*{The Cheng-Yau metric on the cone}
Let $\mathcal{C}\subset {\RR}^3$ be the open cone over the domain $\Omega$ in affine space $E={\RR}^3$, i.e., the cone $\cC=\{(tx,t)\in {\RR}^3|x\in\Omega,t>0\}$. 
Cheng and Yau in \cite{Cheng-Yau 2} show that there exists a unique strictly convex function $\sigma$ on the convex cone $\cC$ satisying  \begin{equation*}
\text{det}(\frac{1}{3}(\log\sigma)_{ij})={\sigma}^2, ~~~\text{and}~~~ \sigma\rightarrow\infty ~~\text{at}~~\partial\Omega,
\end{equation*}
and that the metric $h=\text{Hess}(\frac{1}{3}(\log\sigma))$ on the cone is complete and invariant under linear automorphisms of $\cC$, will be referred as to Cheng-Yau metric.\\
Calabi \cite{Calabi} and Cheng and Yau \cite{Cheng-Yau 1} show that each convex domain $\Omega$ has associated to it a unique strictly convex function $u$ satisfying 
\begin{equation*}
\text{det}(u_{ij})=(\frac{1}{u})^4,~~~\text{and}~~~ u|_{\partial{\Omega}}=0.
\end{equation*}
The radial graph of $-\frac{1}{u}$ is a hyperbolic affine sphere $H$ asymptotic to the boundary of the cone $\mathcal{C}$ with center 0 and affine mean curvature -1. Moreover, the metric $-\frac{1}{u}u_{ij}dt^idt^j$ on $\Omega$ and the associated affine metric $g$ on $H$ are isometric under the map $-\frac{1}{u}(t_1,t_2,1):\Omega\rightarrow H$. (so we may sometimes just write $g=-\frac{1}{u}u_{ij}dt^idt^j$.)
\end{subsection}

The relation between the Cheng-Yau metric $h$ on the cone and the affine metric $g$ on the hyperbolic affine sphere inside the cone is as follows:

\begin{proposition} \label{2}(Loftin \cite{Loftin 3})
Let $h=\frac{1}{3}(d^2\log\sigma)$ be the Cheng-Yau metric on the cone $\cC$, then hypersurface $H={\sigma}^{-1}(1)$ is a hyperbolic affine sphere with center the origin and affine mean curvature -1. The natural foliation $\cC=\cup_{s>0}{sH}$ gives the metric splitting 
\begin{equation*}
(\cC,h)=({\RR}^{+},\frac{ds^2}{s^2})\bigoplus(H,g),
\end{equation*}
where $g$ the Blaschke metric along $H$.
\end{proposition}

By definition, the affine normal (the position vector in the above case) of $H$ is invariant under $Aut(H)$. Combining this with equation (\ref{metric splitting}), the affine metric $g$ on $H$ is also invariant under $Aut(H)$ and hence the metric $g=-\frac{1}{u}u_{ij}dt^idt^j$ on $\Omega$ descends to a metric on $M=\Omega/\Gamma$ which is also called affine or Blaschke metric on $M$.  The Cheng-Yau metric $h$ is invariant under $Aut(\cC)$, and hence is also invariant under $Aut(H)$, since $Aut(H)\subset Aut(\cC)$.

\begin{subsection}*{The construction of the Loftin metric} Loftin mentioned in \cite{Loftin 3} that the construction of the DG metric (will be defined later) can be carried out with other invariant affine K$\ddot{\text{a}}$hler metrics, e.g., Cheng-Yau metric, instead of the Koszul-Vinberg metric on the cone. Hence we define a different Riemannian metric on $\mathcal{B}(S)$ using a construction similar to that of the DG metric but using the Cheng-Yau metric on the cone instead of the Koszul-Vinberg metric. For the above reason, we call the new metric "the Loftin metric".

We begin the construction of the Loftin metric: 

Suppose that the pair $(f,M)\in \mathcal{B}(S)$ corresponds to a convex $\RR P^2$ structure on $S$. Let $\mathcal{C}\subset R^3$ be the corresponding open cone in affine space $E={\RR}^3$. 

On the one hand, the Cheng-Yau metric $h$ at each point of the hyperbolic affine sphere $H$ is an inner product on ${\RR}^3$, hence induces an inner product on $sl(3,\RR)\subset gl(3,\RR)\subset Hom({\RR}^3,{\RR}^3)\cong {\RR}^3\otimes ({\RR}^3)^*$, and therefore also induces a Riemannian metric $l=h\otimes h^*$ on the trivial bundle $sl(3,\RR)\times H$ over $H$. We then obtain a Riemannian metric $l$ (abusing notation) on the bundle $sl(3,\RR)_{Ad \rho}$. Explicitly, supposing that $\phi, {\phi}^{\prime}$ are sections of $sl(3,\RR)_{Ad \rho}$, $\forall p\in S,$ we define \begin{equation}\label{l}l(\phi,\phi^{\prime})|_p:=l_x(\widetilde{\phi}_x,\widetilde{\phi}^{\prime}_x), \text{for some $x\in {\pi}^{-1}(p)\subset H$},\end{equation} where $\pi$ is composition of the projective map from hypersurface $H$ to $\Omega$ and the quotient map from $\Omega$ to $S$, and $\widetilde{\phi},\widetilde{\phi}^{\prime}$ are the liftings of sections $\phi,{\phi}^{\prime}$ of $sl(3,\RR)_{Ad \rho}$ to the trivial bundle $sl(3,\RR)\times H$. Then by definition of the bundle $sl(3,\RR)_{Ad \rho}$ (see equation (\ref{bundle})), we see that (for $\gamma\in \rho(\pi)$, we have that $\widetilde{\phi},\widetilde{\phi}^{\prime}$ satisfy  \begin{equation}\label{equality}\widetilde{\phi}_{\gamma x}=Ad(\gamma)\widetilde{\phi}_{x},~~\text{and}~~ \widetilde{\phi}^{\prime}_{\gamma x}=Ad(\gamma)\widetilde{\phi}^{\prime}_{x}.\end{equation} 

For any $\gamma\in \rho(\pi)=\Gamma<Aut(H)$, we compute \begin{eqnarray*}
&&l_{\gamma x}(\widetilde{\phi}_{\gamma x},\widetilde{\phi}^{\prime}_{\gamma x})\\&=&l_{\gamma x}(Ad(\gamma)\widetilde{\phi}_{x},Ad(\gamma)\widetilde{\phi}^{\prime}_{x})~~\text{by equation (\ref{equality})}\\&&\text {under the identification of $gl(3,\RR)$ and ${\RR}^3\otimes {{\RR}^*}^3$}\\&=&l_{\gamma x}((Ad(\gamma)\widetilde{\phi}_{x})^i_j e_i\otimes e^j, (Ad(\gamma)\widetilde{\phi}^{\prime}_{x})^k_l e_k\otimes e^l)\\
&=&h_{\gamma x}\otimes h^*_{\gamma x}((\widetilde{\phi}_{x})^i_j\gamma e_i\otimes (\gamma e_j)^*,(\widetilde{\phi}^{\prime}_{x})^k_l\gamma e_k\otimes(\gamma e_l)^*)~~~\text{since $l=h\otimes h^*$}\\&=&(\widetilde{\phi}_{x})^i_j(\widetilde{\phi}^{\prime}_{x})^k_l h_{\gamma x}(\gamma e_i, \gamma e_k) h^*_{\gamma x}((\gamma e_j)^*, (\gamma e_l)^*)\\
&& \text{Since $h$ and $h^*$ are invariant under affine automorphisms of $H$,}\\
&& h_{\gamma x}(\gamma e_i, \gamma e_k)=h_{x}(e_i, e_k)~~ \text{and}~~ h^*_{\gamma x}((\gamma e_j)^*, (\gamma e_l)^*)=h^*_{x}(e_j^*, e_l^*)\\
&=&(\widetilde{\phi}_{x})^i_j(\widetilde{\phi}^{\prime}_{x})^k_l h_{x}(e_i, e_k) h^*_{x}(e^j, e^l)~~~\text{noting that}~~ e^j=e_j^*, e^l=e_l^*\\
&=&l_x(\widetilde{\phi}_x,\widetilde{\phi}^{\prime}_x) ~~~\text{since $l=h\otimes h^*$}.
\end{eqnarray*} Hence we obtain that $l(\phi,\phi^{\prime})|_p$ does not depend on the choice of $x$ in $ {\pi}^{-1}(p)$ and hence is well-defined.

On the other hand, the affine metric and the orientation on $S$ define a metric on $\mathcal{A}^p(S)$ (the space of 1-forms on $S$) and hence enable us to define a Hodge star operator 
\begin{equation*}
*:\mathcal{A}^p(S)\rightarrow \mathcal{A}^{2-p}(S)
\end{equation*} by setting
\begin{equation}
\alpha\wedge*\beta=<\alpha,\beta>dvol.\label{Hodge dual}
\end{equation} Combining the action of Hodge star operator with the Riemannian metric $l$ on the bundle $sl(3,\RR)_{Ad \rho}$, we may define a positive definite inner product $\widetilde{g}_{Loftin}$ on the space  $\mathcal{A}^1(S,sl(3,\RR)_{Ad \rho})$, the space of 1-forms taking values in the bundle $sl(3,\RR)_{Ad \rho}$ as follows. 

Let $\sigma\otimes \phi, {\sigma}^{\prime}\otimes{\phi}^{\prime}\in \mathcal{A}^1(S,sl(3,\RR)_{Ad \rho})$, where $\sigma,{\sigma}{\prime}\in \mathcal{A}^1(S)$ and $\phi, {\phi}^{\prime}$ are sections of $ sl(3,\RR)_{Ad \rho}$. We define a pairing $\widetilde{g}_{Loftin}$ as follows:
\begin{equation}\label{widetilde}
\widetilde{g}_{Loftin}(\sigma\otimes \phi, {\sigma}^{\prime}\otimes{\phi}^{\prime})= \int_S(\sigma\wedge*{\sigma}{\prime})l(\phi, {\phi}^{\prime}).
\end{equation} 
By linearity, we may extend the definition of $\widetilde{g}_{Loftin}$ to a pair $\Sigma_i\sigma_i\otimes \phi_i$ and $\Sigma_j{\sigma}_j^{\prime}\otimes{\phi}_j^{\prime}$. Hence we obtain an inner product, which is also denoted $\widetilde{g}_{Loftin}$, defined on the whole space $\mathcal{A}^1(S,sl(3,\RR)_{Ad \rho})$.

Next we define a metric, which will be referred as to the Loftin metric $g_{\text{Loftin}}$, on the cohomology $H^1(S;sl(3,\RR)_{Ad \rho})$ as follows (see pp. 108-111 in M.S.Raghunathan \cite{Raghunathan}): 
The metric $l$ we defined on the fibers of $sl(3,\RR)_{Ad\rho}$ gives an isomorphism 
\begin{eqnarray*}
\sharp:sl(3,\RR)_{Ad\rho}\rightarrow sl(3,{\RR}^*)_{Ad{\rho}^*}, \\
\text{where $({\rho}^*y)(x)=y({\rho}^{-1}x)$, for $y\in {\RR}_3, \text{and} x\in {\RR}^3$}
\end{eqnarray*}
defined by setting 
\begin{equation*}
(\sharp v)_x(u_x)=l_x(u_x,v_x)
\end{equation*}
for $u_x,v_x\in sl(3,\RR)$, and $x\in M$. This isomorphism extends naturally to an isomorphism again denoted $\sharp$ of ${\cA}^p(S,sl(3,\RR)_{Ad\rho})$ on ${\cA}^p(S,sl(3,{\RR}^*)_{Ad{\rho}^*})$:
\begin{equation*}
\sharp:{\cA}^p(S,sl(3,\RR)_{Ad\rho})\rightarrow {\cA}^{p}(S,sl(3,{\RR}^*)_{Ad{\rho}^*}).\label{sharp defn}
\end{equation*}
In addition, the Hodge star operator on ${\cA}^p(S)$ naturally extends to be defined on ${\cA}^p(S,sl(3,\RR)_{Ad\rho})$.
Finally, we define an operator $\delta:{\cA}^p(S,sl(3,\RR)_{Ad\rho})\rightarrow{\cA}^{p-1}(S,sl(3,\RR)_{Ad\rho})$ by setting 
\begin{equation}\label{delta}
\delta=-(\sharp)^{-1}*^{-1}d*(\sharp)
\end{equation}
and then define the Laplacian $\Delta:{\cA}^p(S,sl(3,\RR)_{Ad\rho})\rightarrow {\cA}^{p}(S,sl(3,\RR)_{Ad\rho})$ by setting 
\begin{equation*}
\Delta=d\delta+\delta d.
\end{equation*}
A form $\xi\in {\cA}^p(S,sl(3,\RR)_{Ad\rho})$ is harmonic if $\Delta\xi=0$. 
In particular, if $M$ is compact, the form $\xi$ is harmonic if and only if \begin{equation*}
d\xi=0,~~~~ \delta\xi=0.
\end{equation*}
The kernel ${\mathcal{H}}^{\infty}(S,sl(3,\RR)_{Ad \rho})$ of $\Delta$ and the images of $d:{\mathcal{A}}^0(S,sl(3,\RR)_{Ad \rho})\rightarrow{\mathcal{A}}^1(S,sl(3,\RR)_{Ad \rho})$
and
$\delta:{\mathcal{A}}^2(S,sl(3,\RR)_{Ad \rho})\rightarrow{\mathcal{A}}^1(S,sl(3,\RR)_{Ad \rho})$ 
decompose the vector space of 1-forms valued in $sl(3,\RR)_{Ad \rho}$ into an orthogonal direct sum
\begin{equation*}
{\mathcal{A}}^1(S,sl(3,\RR)_{Ad \rho})={\mathcal{H}}^{\infty}(S,sl(3,\RR)_{Ad \rho})\oplus
Image(d)\oplus Image(\delta)
\end{equation*}

Since each de Rham cohomology class contains a unique harmonic representative from non-abelian Hodge theory (see Proposition 7.10 in \cite{Raghunathan}), we may define the pairing at $[\rho]$
\begin{equation*}
g_{\text{Loftin}}: {H}^{1}(S,sl(3,\RR)_{Ad \rho})\times{H}
^{1}(S,sl(3,\RR)_{Ad \rho})\rightarrow R
\end{equation*}
\begin{equation*}
\text{by } g_{\text{Loftin}}([\alpha], [\beta]):=\widetilde{g}(\alpha_{harm},\beta_{harm}), ~~\text{for} [\alpha], [\beta]\in {H}^{1}(S,sl(3,\RR)_{Ad \rho}),
\end{equation*} 
where $\alpha_{harm},\beta_{harm}$ are the unique harmonic representatives of $[\alpha], [\beta]$ respectively and $\widetilde{g}$ is defined above (see equation (\ref{widetilde})). Hence we have a well-defined Riemannian metric $g_{\text{Loftin}}$ on Goldman space $\mathcal{B}(S)$.
\end{subsection}\\

Next we introduce the Riemannian metric on the space $\mathcal{B}(S)$ defined by Darvishzadeh and Goldman (see details in \cite{Goldman 2}). Since this Riemannian metric relies heavily on the Koszul-Vinberg metric on the cone $\cC$ in ${\RR}^3$, we first recall the definition of the Koszul-Vinberg metric.

\begin{subsection}*{Koszul-Vinberg metric}
Let $\mathcal{C}\subset {\RR}^3$ be a cone in affine space $E={\RR}^3$. The dual cone ${\mathcal{C}}^*$ is the subset of the dual vector space $E^*$ consisting of linear functionals $\psi:E\rightarrow \RR$ which are positive on $\mathcal{C}$.

Recall the Koszul-Vinberg characteristic function $k(x)$ on the cone $\mathcal{C}$: for $x\in \mathcal{C}$, define 
\begin{equation*}
k(x)=\int_{\mathcal{C}^*}e^{-\psi(x)}d\psi.
\end{equation*}
Note that 
\begin{eqnarray}
k(\gamma x)&=&\int_{\mathcal{C}^*}e^{-\psi(\gamma x)}d\psi\nonumber\\
&=&\int_{\mathcal{C}^*}e^{-{\gamma}^*\psi(x)}d\psi\nonumber\\
&=&\int_{\mathcal{C}^*}e^{-\phi(x)}det(\gamma)^{-1}d\phi, \text{~~let $\phi={\gamma}^*\psi\in \mathcal{C}^*, d\phi=det(\gamma)d\psi$}\nonumber\\
&=&det(\gamma)^{-1}\int_{\mathcal{C}^*}e^{-\phi(x)}d\phi\nonumber\\
&=&det(\gamma)^{-1}k(x)\text{, for any $\gamma\in$ Aut($\mathcal{C}$).}\label{k}
\end{eqnarray} Hence the Hessian $d^2\log k$ is  invariant under Aut($\mathcal{C}$). Moreover, the Hessian $d^2\log k$ is actually a positive definite symmetric bilinear form $h$ (which we call the Koszul-Vinberg metric) on the cone $\cC$. 
\end{subsection}

\begin{subsection}*{The construction of DG metric}
We assume the same notation $M, \Omega,  \Gamma,$ and $\cC$ as in the construction of the Loftin metric.

On one hand, for every $x\in k^{-1}(1)$, the Koszul-Vinberg metric at each point of the cone gives an inner product on ${\RR}^3$, hence induces an inner product on $sl(3,\RR)\subset gl(3,\RR)\cong Hom(\RR^3,\RR^3)\cong{\RR}^3\otimes {{\RR}^{*}}^3$, and therefore also induces a Riemannian metric on the bundle $sl(3,\RR)_{\text{Ad}\rho}$, since the Koszul-Vinberg metric is invariant under Aut($\mathcal{C}$).

On the other hand, consider the map $m:\Omega\rightarrow \mathcal{C}$, which takes $[p]\longmapsto k(p)^{\frac{1}{3}}p$. Setting $t$ a positive constant, we have that $tI$ is an element of Aut($\mathcal{C}$). After substituting $tI$ for $\gamma$ into equation (\ref{k}), we obtain that \begin{equation}
k(tp)=t^{-3}k(p).\label{ke}
\end{equation} Then $k(tp)^{\frac{1}{3}}tp=k(p)^{\frac{1}{3}}p$ and hence $m$ is well-defined (i.e., $k(p)^{\frac{1}{3}}p$ does not depend on the choice of elements in $[p]$). Moreover, after substituting $k(p)^{\frac{1}{3}}$ for $t$ into equation (\ref{ke}), we obtain that $k(k(p)^{\frac{1}{3}}p)=k(p)^{-1}k(p)=1$ and hence $m(\Omega)=k^{-1}(1)$. 

The Riemannian metric $m^*(d^2\log k)$ on $\Omega$ is invariant under $\Gamma$. Hence $m^*(d^2\log k)$ defines a Riemannian metric on $\Omega/{\Gamma}$. Thus corresponding to every convex $\RR P^2$-structure on $S$, there exists an associated Riemannian metric on $S$. Now the remainder is similar to the definition of the Loftin metric: we first have the induced metric $\widetilde{g}_{DG}$ on the space $\mathcal{A}^1(S;sl(3,\RR)_{Ad \rho})$ and then define Laplacian operator on the space, hence obtain harmonic representatives as kernel of the Laplacian operator. The DG metric $g_{DG}$ on the cohomologous classes as tangent vectors is actually defined as the metric $\widetilde{g}_{DG}$ on harmonic representatives.
\end{subsection}
\end{section}
\begin{section}{Embedding of Teichm$\ddot{\text{u}}$ller Space}\label{e}
Inside this section, we fix the notation as follows:\\
(i) the hyperbolic affine sphere( the hyperboloid) $H=\{{x_3}^2-{x_1}^2-{x_2}^2=1\}$;\\
(ii) the domain $\Omega=\{{t_1}^2+{t_2}^2<1\}\subset {\RR}^2$.

By definition, Goldman space $\cB(S)$ is the space of all convex real projective structures on the surface and we can think of the Teichm$\ddot{\text{u}}$ller locus inside $\cB(S)$ as the subspace of convex real projective structures which arise from hyperbolic structures. Noting that there are a variety of viewpoints of the deformation space $\cB(S)$, in this section we give a detailed description of the embedding of Teichm$\ddot{\text{u}}$ller space $\cT(S)$ from some different viewpoints and then show their equivalence.
\begin{enumerate}
\item When Goldman space $\cB(S)$ is identified with the space of affine sphere structures that can be given on the surface $S$ (see Proposition \ref{affine}), then the Teichm$\ddot{\text{u}}$ller locus consists of points representing surfaces which admits an affine sphere structure whose affine sphere is the hyperboloid $H$;\\
\item When Goldman space $\cB(S)$ is identified with the space containing pairs $(\Sigma,U)$ (see Proposition \ref{U}), then the Teichm$\ddot{\text{u}}$ller locus is the subspace containing pairs $(\Sigma,0)$;\\
\item When Goldman space $\cB(S)$ is identified with an open subspace of the space containing conjugate classes of representations $\rho: \pi\rightarrow SL(3,\RR)$, then the Teichm$\ddot{\text{u}}$ller locus in $\mathcal{B}(S)$ consists exactly of conjugation classes of representations ${\rho}^{\prime}:\pi\rightarrow PSL(2,\RR)$ after composing by the irreducible representation $\Phi: PSL(2,\RR)\rightarrow SL(3,\RR)$.
\end{enumerate}

Because of the equivalence of the definitions described above, in following sections we will use the description of the embedding of Teichm$\ddot{\text{u}}$ller space from different viewpoints for convenience without explanation.\\

Firstly note that (2) follows from (1) immediately, because the Pick form $C$ for the hyperboloid vanishes, hence the cubic differential is 0 (see the end of \S \ref{correspondence}). Then we continue to describe (1) and (3).\\

To explain (1): 
Consider a hyperbolic structure on the surface $S$, then the hyperbolic surface $S$ has ${\HH}^2$ (the upper half plane in $\CC$ with hyperbolic metric $\frac{1}{y^2}|dz|^2$) as its Riemannian cover. The following lemma shows that the hyperbolic metric can be realized as the Blaschke metric on $\Omega$ induced from the hyperbolic affine sphere $H$. Hence the hyperbolic surface $S$ actually admits an affine sphere structure as a quotient of the hyperbolic affine sphere $H$.
\begin{lemma}(Kim \cite{Kim})\label{embedding} Suppose the domain $\Omega$ is given with Blaschke metric and ${\HH}^2$ is given with the hyperbolic metric, then the map defined by
\begin{equation*}
F:\Omega\rightarrow {\HH}^2
\end{equation*}
\begin{equation*}
F(t_1,t_2)=\frac{t_1}{1-t_2}+i\frac{1}{1-t_2}\sqrt{1-{t_1}^2-{t_2}^2}
\end{equation*}~~~~~~~~~  is an isometry.
\end{lemma}
\textit{Proof}.
Since $x=\frac{t_1}{1-t_2}, ~y=\frac{1}{1-t_2}\sqrt{1-{t_1}^2-{t_2}^2}$, then 
\begin{eqnarray}
dx&=&\frac{1}{1-t_2}dt_1+\frac{t_1}{{1-t_2}^2}dt_2, \label{dx} \\
dy&=&\frac{-t_1}{(1-t_2)\sqrt{1-{t_1}^2-{t_2}^2}}dt_1+\frac{1-{t_1}^2-t_2}{(1-t_2)^2\sqrt{1-{t_1}^2-{t_2}^2}}dt_2 \label{dy}.
\end{eqnarray}
We note that the function $u=-\sqrt{1-{t_1}^2-{t_2}^2}$ on $\Omega$ is the solution to the equation
\begin{equation*}
\text{det}(\frac{{\partial}^2 u}{\partial t^i\partial t^j})=(\frac{1}{u})^4.
\end{equation*}
Then by applying equations (\ref{dx}) and (\ref{dy}), we compute the hyperbolic metric on ${\HH}^2$:
\begin{eqnarray*}
&&\frac{1}{y^2}(dx^2+dy^2)\\ 
&=&\frac{1-{t_2}^2}{(1-{t_1}^2-{t_2}^2)^2}{dt_1}^2+\frac{2t_1 t_2}{(1-{t_1}^2-{t_2}^2)^2}{dt_1dt_2}+\frac{1-{t_1}^2}{(1-{t_1}^2-{t_2}^2)^2}{dt_2}^2\\
&&\text{From the solution $u(t_1,t_2)=-\sqrt{1-{t_1}^2-{t_2}^2}$},\\
&=& -\frac{1}{u}u_{11}{dt^1}^2 -2\frac{1}{u}u_{12}{dt^1}{dt^2}-\frac{1}{u}u_{11}{dt^2}^2,\\
&=& -\frac{1}{u}u_{ij}dt^idt^j, ~~\text{the Blaschke metric on}~ \Omega. \qed
\end{eqnarray*}
\textit{Remark.} Instead of $(\Omega, \text{the Blaschke metric})$, Kim \cite{Kim} actually defined the map $F$ in the lemma on $(\Omega, \text{the Hilbert metric})$. But in fact the Blaschke metric and the Hilbert metric are the same in this case.  Hence our approach in the lemma above is a bit different from his proof in \cite{Kim}. In the end, we can compute 
\begin{equation} F^{-1}(x,y)= (\frac{2x}{x^2+y^2+1}, 1-\frac{2}{x^2+y^2+1}).\label{F^{-1}}
\end{equation}\\
In the following lemma, we give a conformal parametrisation of the hyperbolic affine sphere $H$ on ${\HH}^2$ (see [\ref{parametrization}]).
\begin{lemma} \label{f} The map $f:{\HH}^2\rightarrow H\subset {\RR}^3$ defined as 
\begin{equation}
f(z)=(\frac{x}{y}, \frac{x^2+y^2-1}{2y},\frac{x^2+y^2+1}{2y})\label{f(z)}
\end{equation} is an isometry.
\end{lemma}
\textit{Proof.}  Once again, set $u=u(t_1,u_2)=-\sqrt{1-{t_1}^2-{t_2}^2}$. Noting that the hypersurface $H$ is exactly the radial graph of the function $-\frac{1}{u}$ (i.e., the image of $-\frac{1}{u}(t_1,t_2,1))$, we have the map $G:(\Omega, \text{the Hilbert metric}) \rightarrow (H, \text{the Blaschke metric})$ is an isometry map, where \begin{equation}
G(t_1,t_2):=-\frac{1}{u}(t_1,t_2,1)=\frac{-1}{\sqrt{1-{t_1}^2-{t_2}^2}}(t_1,t_2,1).\label{U}
\end{equation} Hence, combining with Lemma \ref{embedding}, we obtain that the composition map $f:=G\circ F^{-1}$ is an isometry from ${\HH}^2\rightarrow H$. Explicitly, we have \begin{eqnarray*}
&&f(z)=G\circ F^{-1}(z)\\
&=&G(\frac{2x}{x^2+y^2+1}, 1-\frac{2}{x^2+y^2+1})~~\text{by definition (equation (\ref{F^{-1}})) of $F^{-1}$,}\\
&=&(\frac{x}{y}, \frac{x^2+y^2-1}{2y},\frac{x^2+y^2+1}{2y})~~\text{by the definition (equation (\ref{U})) of $G$.}
\end{eqnarray*} \qed

To explain (3): A hyperbolic structure on $S$ determines a holonomy homomorphism $\pi\rightarrow PSL(2,\RR)$. Elements in $PSL(2,\RR)$ keeps the hyperbolic metric invariant. So we hope to find $\Phi: PSL(2,\RR)\rightarrow SL(3,\RR)$ such that image of $\Phi$ fixes the hyperbolic affine sphere $H$ and the Blaschke metric along $H$. Equivalently, we wish to show that $f$ in Lemma \ref{f} is $\Phi$-invariant, since the map $f$ is an isometry from ${\HH}^2$ to $H\subset {\RR}^3$. We eventually realize the hyperbolic structure on $S$ as a convex real projective structure with the map $\Phi$ defined in the following Proposition:
\begin{proposition}\label{inside}(Kim \cite{Kim}) The map $\Phi: PSL(2,\RR)\rightarrow SL(3,\RR)$ defined by \\
\[ \Phi(A)=\left( \begin{array}{ccc}
ad+bc & ac-bd & ac+bd \\
ab-cd & \frac{a^2-b^2-c^2+d^2}{2}& \frac{a^2+b^2-c^2-d^2}{2} \\
ab+cd & \frac{a^2-b^2+c^2-d^2}{2}&\frac{a^2+b^2+c^2+d^2}{2} \end{array} \right)\] for \[A=\left( \begin{array}{ccc}
a&b\\c&d\end{array} \right)\]\\
is an injective homomorphism of $PSL(2,\RR)$ into $SL(3,\RR)$ with image $SO(2,1)$ such that the map $f$ in Lemma \ref{f} is $\Phi$-equivariant. 
\end{proposition}

\textit{Remark.}    \label{Remark}
The map $\Phi$ induces a Lie algebra homomorphism at the identity matrix. Abusing the notation $\Phi$, we obtain that
\[ \Phi(A)=\left( \begin{array}{ccc}
0 & c-b & c+b \\
b-c & 0& 2a\\
b+c & 2a&0\end{array} \right)\] for \[A=\left( \begin{array}{ccc}
a&b\\c&-a\end{array} \right)\in sl(2,\RR),\]\\
which is a Lie algebra homomorphism of $sl(2,\RR)$ into $sl(3,\RR)$ with image $so(2,1)\subset sl(3,\RR)$. This map $\Phi$ on the Lie algebra will help us connect tangent vectors of Teichm$\ddot{\text{u}}$ller space with tangent vectors of the Teichm$\ddot{\text{u}}$ller locus in $\cB(S)$.
\end{section}
\begin{section}{The restriction of two generalized Weil-Petersson metrics}\label{restriction}
The goal of this section is to prove 
\begin{theorem}\label{WP} The DG metric and the Loftin metric both restrict to a constant multiple of the Weil-Petersson metric on Teichm$\ddot{\text{u}}$ller space.
\end{theorem}
It is clear that Theorem \ref{WP}  follows from the following two propositions. 

\begin{proposition}\label{main} The restriction of the DG metric $g_{DG}$ to Teichm$\ddot{\text{u}}$ller space is a constant multiple of the restriction of the Loftin metric $g_{Loftin}$ to Teichm$\ddot{\text{u}}$ller space.
\end{proposition}
\begin{proposition}\label{Loftin metric} The Loftin metric $g_{Loftin}$ restricts to a constant multiple of the Weil-Petersson metric on Teichm$\ddot{\text{u}}$ller space.
\end{proposition}
Hence the remaining goal of this section is to show the above two propositions. We first finish the proof of Proposition \ref{main} and then show Proposition \ref{Loftin metric}.

\begin{subsection}*{Proof of Proposition \ref{main}} 
We start with showing Lemma \ref{1}, which is essential in the proof of Proposition \ref{main} which compares the Cheng-Yau metric with the Koszul-Vinberg metric on the cone $\cC=\{{x_3}^2>{x_1}^2+{x_2}^2\}$.

\begin{lemma}(Sataki \cite{Sataki})\label{1}  In the case of the cone $\mathcal{C}=\{{x_3}^2>{x_1}^2+{x_2}^2\}$, the function $k=({x_3}^2-{x_1}^2-{x_2}^2)^{-\frac{3}{2}}$ is the characteristic function of the cone $\cC$. Supposing the function $\sigma$ is the solution to the Cheng-Yau equation \begin{equation*}
\text{det}(\frac{1}{3}(\log\sigma)_{ij})={\sigma}^2 \text{~~~~~on the cone $\mathcal{C}$},
\end{equation*} then \\
(i) we have that $k={\sigma}$;\\
(ii) the Koszul-Vinberg metric on the cone is 3 times the Cheng-Yau metric; and  \\
(iii) the hypersurface $k^{-1}(1)$ coincides with the hypersurface $\sigma^{-1}(1)$, the hyperbolic affine sphere which is asymptotic to the boundary of the cone. \\
(iv) the metric on $M=\Omega/\Gamma$ where $\Omega=\{{t_1}^2+{t_2}^2<1\}\subset\subset {\RR}^2\subset {\RR}P^2$,  induced from Koszul-Vinberg metric on the cone,  is 3 times the affine metric obtained from immersing the hyperbolic affine sphere $\sigma^{-1}(1)$.
\end{lemma}
Remark. Part (i), (ii), (iii) of this lemma is already proved in by Sataki in \cite{Sataki} (the statement is a bit different from the original version), here we give a detailed computation to prove it.

\textit{Proof}. 
(i): In the cone $\mathcal{C}=\{{x_3}^2>{x_1}^2+{x_2}^2\}$,  the characteristic function is $k(x)=({x_3}^2-{x_1}^2-{x_2}^2)^{-\frac{3}{2}}$ (see Example 4.2, pp 67 in \cite{Shima}). \\
Then we have the following computation 
\begin{eqnarray}
&&\text{det} (\frac{1}{3}\frac{{\partial}^2 \log k}{\partial {x^i}\partial {x^j}})\nonumber\\
&=&\text{det} (\frac{1}{3}\frac{{\partial}^2 \log({x_3}^2-{x_1}^2-{x_2}^2)^{-\frac{3}{2}}}{\partial {x^i}\partial {x^j}})\nonumber\\
&=&\text{det} (-\frac{1}{2}\frac{{\partial}^2 \log({x_3}^2-{x_1}^2-{x_2}^2)}{\partial {x^i}\partial {x^j}}) \nonumber\\
&=&\text{det} \begin{pmatrix}\frac{2{x_1}^2+t}{t^2}&\frac{2{x_1}{x_2}}{t^2}&\frac{-2{x_1}{x_3}}{t^2}\\\frac{2{x_1}{x_2}}{t^2}&\frac{2{x_2}^2+t}{t^2}&\frac{-2{x_2}{x_3}}{t^2}\\\frac{-2{x_1}{x_3}}{t^2}&\frac{-2{x_2}{x_3}}{t^2}&\frac{2{x_3}^2-t}{t^2}\end{pmatrix}\label{hh}\\
&=&t^{-3},\text{after setting $t={x_3}^2-{x_1}^2-{x_2}^2$}\nonumber\\
&=&k^2, \text{noting that $k=t^{-\frac{3}{2}}$}.\nonumber
\end{eqnarray}
Therefore by the uniqueness of solution of the Cheng-Yau equation, the function $k(x)=({x_3}^2-{x_1}^2-{x_2}^2)^{-\frac{3}{2}}$ coincides with the solution $\sigma$ to the Cheng-Yau equation  \begin{equation*}
\text{det}(\frac{1}{3}(\log\sigma)_{ij})={\sigma}^2.
\end{equation*} This finishes the proof of (i).\\
(ii): From the definition of the Koszul-Vinberg metric as $d^2(\log k)$ and the Cheng-Yau metric as $\frac{1}{3}(d^2\log\sigma)$, and combining with the equality $k=\sigma$ in (i), we conclude the proof of (ii). \\
(iii) immediately follows from the equality $k=\sigma$ in (i).\\
(iv): Combine the fact that the metric on $M=\Omega/\Gamma$ is the quotient metric of the Koszul-Vinberg metric restricted to the hypersurface $k^{-1}(-1)$ and the affine metric on $M$ is the quotient metric of the Cheng-Yau metric restricted to the hypersurface ${\sigma}^{-1}(-1)$ with the the facts in (ii) and (iii), we conclude the proof of (iv).  \qed\\

For further reference, we collect the concepts and facts associated to the convex ${\RR P}^2$-structure $M$ on the surface arising from the hyperbolic structures. \\
(\MakeUppercase{\romannumeral 1}) the domain $\Omega=\{{t_1}^2+{t_2}^2<1\}\subset {\RR}^2\subset {\RR P}^2$;\\
\MakeUppercase{(\romannumeral 2)} the holonomy group $\Gamma<Aut(\cC)<SL(3,\RR)$;\\
\MakeUppercase{(\romannumeral 3)} the manifold $M\cong \Omega/{\Gamma}$;\\
\MakeUppercase{(\romannumeral 4)} the solution $\sigma(x)=({x_3}^2-{x_1}^2-{x_2}^2)^{-\frac{3}{2}}$ to the Cheng-Yau equation on the cone $\cC$ \begin{equation*}
\text{det}(\frac{1}{3}(\log\sigma)_{ij})={\sigma}^2;
\end{equation*}
\MakeUppercase{(\romannumeral 5)} the solution $u(t_1,t_2)=-\sqrt{1-{t_1}^2-{t_2}^2}$ to the Cheng-Yau equation on the domain $\Omega$ \begin{equation*}
\text{det}(\frac{{\partial}^2 u}{\partial t^i\partial t^j})=(\frac{1}{u})^4;
\end{equation*}
\MakeUppercase{(\romannumeral 6)} the hyperbolic affine sphere $H$ asymptotic to the boundary of the cone $\mathcal{C}$ with center 0 and constant affine mean curvature -1 is exactly the hypersurface ${\sigma}^{-1}(1)$, and also the radial graph of the function $-\frac{1}{u}$ (i.e., the image of $-\frac{1}{u}(t_1,t_2,1)$).\\

\textit{Proof of Proposition \ref{main}}. We start by noting that the Teichm$\ddot{\text{u}}$ller locus in $\mathcal{B}(S)$ exactly contains the convex $\RR P^2$ manifolds which are diffeomorphic to $\Omega/{\Gamma}$, where $\Omega=\{{t_1}^2+{t_2}^2<1\}$. 

Then by comparing the two different pairings $\widetilde{g}_{DG}$ and  $\widetilde{g}_{Loftin}$ on the space ${\cA}^1(S,sl(3,\RR)_{Ad \rho})$, and by (ii),(iii) and (iv) of Lemma \ref{1}, we obtain that they are isometric when restricted to Teichm$\ddot{\text{u}}$ller locus (up to a constant mulitiple). 

Next, since the pairing $\widetilde{g}_{DG}$ and the pairing $\widetilde{g}_{Loftin}$ are isometric (up to a constant multiple), then the harmonic representatives ${\alpha}^{DG}_{harm}$ and ${\alpha}^{Loftin}_{harm}$ in the cohomology class $[\alpha]\in {H}^1(S;sl(3,\RR)_{Ad \rho})$, in the kernel of Laplacian operators for the metrics $\widetilde{g}_{DG}$ and $\widetilde{g}_{Loftin}$ respectively on the space ${\cA}^1(S,sl(3,\RR)_{Ad \rho})$, are the same.

Finally, recall that \begin{eqnarray*}
&g_{DG}([\alpha], [\beta])&=\widetilde{g}_{DG}({\alpha}^{DG}_{harm},
{\beta}^{DG}_{harm})\\ \text{and}& g_{Loftin}([\alpha], [\beta])&=\widetilde{g}_{Loftin}({\alpha}^{Loftin}_{harm}, {\beta}^{Loftin}_{harm}), 
\end{eqnarray*} for $[\alpha], [\beta]\in {H}^1(S;sl(3,\RR)_{Ad \rho})$. From above arguments, we then conclude the pairing $g_{DG}([\alpha],[\beta])$ and $g_{Loftin}([\alpha],[\beta])$ give the same result (up to a constant multiple).\qed
\end{subsection}
\begin{subsection}*{Proof of Proposition \ref{Loftin metric}}
We want to show that the Loftin metric restricts to be the Weil-Petersson metric on Teichm$\ddot{\text{u}}$ller space. It requires some preparation to achieve this goal. We need to understand the following two objects:
\begin{enumerate}
\item an explicit description of the metric $l$ on the Lie algebra bundle $sl(3,\RR)_{Ad \rho}$.
\item harmonic representatives in the cohomology class which are tangent vectors on Teichm$\ddot{\text{u}}$ller space, since we define the Loftin metric after choosing the harmonic representative in its cohomology class. 
\end{enumerate}
\begin{subsection}*{(1) The Metric $l$ on the Lie Algebra Bundle}
The main goal of this part is to give an explicit formula for the metric $l$ on the bundle $sl(3,\RR)_{Ad{\rho}}$ by Lemma \ref{pairing g} and hence to compute the norm of $\Phi\bigl(\begin{smallmatrix}-z&z^2\\ -1&z\end{smallmatrix}\bigr)$ ($\Phi$ is the Lie algebra homomorphism of $sl(2,\RR)$ into $sl(3,\RR)$ defined in the remark in the end of \S \ref{e}) in Corollary \ref{pairing}, which is useful in step (iii) of the proof of Proposition \ref{Loftin metric}. 
\begin{lemma}\label{pairing g} Supposing that the metric $l$ is defined on the bundle $sl(3,\RR)_{Ad \rho}$, the matrix $h$ is the matrix presentation of the Cheng-Yau metric at point $p$ of the hyperbolic affine sphere $H$ under the standard basis of ${\RR}^3$(i.e., $e_1=(1,0,0)^T, e_2=(0,1,0)^T,e_3= (0,0,1)^T$, and $h(e_i,e_j)=h_{ij}$), we then obtain that  \begin{equation*}
l_{p}(A,B)=tr(A^Th^{-1}Bh) ~~~~~~~~~~~\text{for $A, B\in sl(3,\RR)$}.
\end{equation*}
\end{lemma}
\textit{Proof}. Recall that $l=h\otimes h^*$. Suppose we have the matrix presentation of the Cheng-Yau metric $h$ in the standard basis as $(h_{ij})$ and the dual basis for $\{e_1,e_2,e_3\}$ is $\{e^1,e^2,e^3\}$ and the matrix presentation of the inverse Cheng-Yau metric $h^{-1}$ is $(h^{ij})$ with $h^{ij}=h^{-1}(e^i,e^j)$. 

Next assuming that the matrix $A=(a_i^j), B=(b_k^l)$ ($i,k$ denote for the row), we obtain that $Ae_i=\mathop{\sum_{\{j=1,2,3\}}}(a_i^j)e_j, Be_k=\mathop{\sum_{\{l=1,2,3\}}}(b_k^l)e_l$. Then we may identify $A$ with $\mathop{\sum_{\{i,j=1,2,3\}}}a_i^j e_j \otimes e^i$ and $B$ with $\mathop{\sum_{\{i,j=1,2,3\}}}{b_k^l e_l \otimes e^k}$, which is exacly the identification of $Hom({\RR}^3,{\RR}^3)\cong {\RR}^3\otimes {{\RR}^*}^3$.
\begin{eqnarray*}
\text{Hence}~~~&&l_{p}(A,B)\\
&=&l_{p}(\mathop{\sum_{\{i,j=1,2,3\}}}a_i^j e_j \otimes e^i, \mathop{\sum_{\{k,l=1,2,3\}}}b_k^l  e_l \otimes e^k)\\
&=&\mathop{\sum_{\{i,j,k,l=1,2,3\}}}a_i^jb_k^ll(e_j\otimes e^i,e_l\otimes e^k)~~\text{by the linearility of $l$}\\
&=&\mathop{\sum_{\{i,j,k,l=1,2,3\}}}a_i^jb_k^lh(e_j,e_l)h^{-1}(e^i,e^k)~~\\
&&(\text{since $l=h\otimes h^*$}), \\
&=&\mathop{\sum_{\{i,j,k,l=1,2,3\}}}a_i^jb_k^lh_{jl}h^{ik}\\
&=&\mathop{\sum_{\{i,j,k,l=1,2,3\}}}a_i^jh^{ik}b_k^lh_{jl}\\
&=&tr(A^Th^{-1}Bh)~~\text{from $A=(a_i^j), B=(b_k^l)$}.\qed
\end{eqnarray*}

\begin{corollary} \label{pairing} Let $\mathcal{C}=\{{x_3}^2>{x_1}^2+{x_2}^2\}$, and $p=f(z)\in H$ for some $z=x+iy\in \HH$, after extending the definition of $l$ in Lemma \ref{pairing g} by  $l_{p}(A, B)=tr(A^Th^{-1}{\overline{B}}h)$, for $A,B\in sl(3,\CC)$,
we obtain \begin{equation*}l_{p}(\Phi\bigl(\begin{smallmatrix}-z&z^2\\ -1&z\end{smallmatrix}\bigr),\Phi\bigl(\begin{smallmatrix}-z&z^2\\ -1&z\end{smallmatrix}\bigr))=16y^{2}.
\end{equation*}
\end{corollary}
\textit{Proof}. Lemma \ref{inside} implies that 
\begin{equation*}
\Phi\begin{pmatrix}-z&z^2\\ -1&z\end{pmatrix}=\begin{pmatrix}0&-z^2-1&z^2-1\\ z^2+1&0&-2z\\z^2-1&-2z&0\end{pmatrix}, \text{denoted as $A$.}
\end{equation*} 

To calculate $l_{f(z)}(A, A)=tr(A^Th^{-1}{\overline{A}}h)$, we need know the matrix presentation of the Cheng-Yau metric $h$, and also $h^{-1}$.

We are in the case of the cone $\mathcal{C}=\{{x_3}^2>{x_1}^2+{x_2}^2\}$, hence the Cheng-Yau metric $h=\frac{1}{3}d^2\log\sigma$, for $\sigma=({x_3}^2-{x_1}^2-{x_2}^2)^{-\frac{3}{2}}$(see Lemma \ref{1}). Since we now restrict ourselves to the hypersurface ${\sigma}^{-1}(1)$, i.e. $\{{x_3}^2-{x_1}^2-{x_2}^2=1\}$, hence the function $t$ in the equation (\ref{hh}) is identically 1, and then we have an explicit formula for the Cheng-Yau metric $h$ and hence $h^{-1}$ as follows, \\\\
$h=\begin{pmatrix}2{x_1}^2+1&2{x_1}{x_2}&-2{x_1}{x_3}\\2{x_1}{x_2}&2{x_2}^2+1&-2{x_2}{x_3}\\-2{x_1}{x_3}&-2{x_2}{x_3}&2{x_3}^2-1\end{pmatrix}$,
$h^{-1}=\begin{pmatrix}2{x_1}^2+1&2{x_1}{x_2}&2{x_1}{x_3}\\2{x_1}{x_2}&2{x_2}^2+1&2{x_2}{x_3}\\2{x_1}{x_3}&2{x_2}{x_3}&2{x_3}^2-1\end{pmatrix}$.\\
Applying the equation (\ref{f(z)}), i.e., $f(z)=(x_1,x_2,x_3)=(\frac{x}{y}, \frac{x^2+y^2-1}{2y},\frac{x^2+y^2+1}{2y})$ to the above matrices, we obtain \\\\
$h=\begin{pmatrix}\frac{2x^2}{y^2}+1&\frac{x(x^2+y^2-1)}{y^2}&-\frac{x(x^2+y^2+1)}{y^2}\\\frac{x(x^2+y^2-1)}{y^2}&\frac{(x^2+y^2-1)^2}{2y^2}+1&-\frac{(x^2+y^2-1)(x^2+y^2+1)}{2y^2}\\-\frac{x(x^2+y^2+1)}{y^2}&-\frac{(x^2+y^2-1)(x^2+y^2+1)}{2y^2}&\frac{(x^2+y^2+1)^2}{2y^2}-1\end{pmatrix}$ and \\\\
$h^{-1}=\begin{pmatrix}\frac{2x^2}{y^2}+1&\frac{x(x^2+y^2-1)}{y^2}&\frac{x(x^2+y^2+1)}{y^2}\\\frac{x(x^2+y^2-1)}{y^2}&\frac{(x^2+y^2-1)^2}{2y^2}+1&\frac{(x^2+y^2-1)(x^2+y^2+1)}{2y^2}\\\frac{x(x^2+y^2+1)}{y^2}&\frac{(x^2+y^2-1)(x^2+y^2+1)}{2y^2}&\frac{(x^2+y^2+1)^2}{2y^2}-1\end{pmatrix}$.\\

Then finally we compute 
\begin{eqnarray*}
l_{f(z)}(A,A)&=&tr(A^Th^{-1}{\overline{A}}h)\\
&=&16y^{2}\text{, where $z=x+iy\in \HH$.}\qed
\end{eqnarray*}

\begin{subsection}*{(2) Harmonic Representative}
The goal of this part is to show 
\begin{lemma}\label{harmonic} Given a Riemann surface $\Sigma$ and $\Phi$: $sl(2,\RR)\rightarrow sl(3,\RR)$ defined in the Remark of Lemma \ref{inside}, then\\
 (i) the tangent space at the point $(\Sigma,0)$ of the image of Teichm$\ddot{\text{u}}$ller space in $\cB(S)$ is exactly spanned by the cohomology class of 
$\phi(z)dz\otimes\Phi\bigl(\begin{smallmatrix}-z&z^2\\-1&z\end{smallmatrix} \bigr)$, where $\phi(z)dz^2$ is a holomorphic quadratic differential on $\Sigma$; and\\
 (ii) the $sl(3,\RR)_{Ad\rho}$-valued 1-forms of the form $\phi(z)dz\otimes\Phi\bigl(\begin{smallmatrix}-z&z^2\\-1&z\end{smallmatrix} \bigr)$ are harmonic representatives (in the sense of the Loftin metric on $ {\mathcal{A}}^{1}(S;sl(3,\RR)_{Ad \rho})$) in their own cohomology class. 
\end{lemma}
\textit{Proof}. (i): Firstly, we note that the image of Teichm$\ddot{\text{u}}$ller space  in $\cB(S)$ exactly contains the representations $\rho:\pi\rightarrow SL(3,\RR)$ which are a composition of ${\rho}^{\prime}:\pi\rightarrow PSL(2,\RR)$ with $\Phi:PSL(2,\RR)\rightarrow SO(2,1)\subset SL(3,\RR)$. Combining this with the fact that the tangent space of Teichm$\ddot{\text{u}}$ller space at $\Sigma$ contains exactly the cohomology classes of $sl(2,\RR)_{Ad{\rho}^{\prime}}$-valued 1-forms of the form $\phi(z)dz\otimes\bigl(\begin{smallmatrix}-z&z^2\\-1&z\end{smallmatrix} \bigr)$ (see \cite{Goldman 1} for details), we conclude the statement of (i).

(ii): (a) We first note that $\phi(z)dz\otimes\Phi\bigl(\begin{smallmatrix}-z&z^2\\-1&z\end{smallmatrix} \bigr)\in {\mathcal{A}}^{1}(S;sl(3,\RR)_{Ad \rho})$ is closed. This follows from the computation
\begin{eqnarray*}
&&d[\phi(z)dz\otimes\Phi\bigl(\begin{smallmatrix}-z&z^2\\-1&z\end{smallmatrix} \bigr)]\\
&=&d[z^2\phi(z)dz\otimes\Phi\bigl(\begin{smallmatrix}0&1\\0&0\end{smallmatrix} \bigr)-\phi(z)dz\otimes\Phi\bigl(\begin{smallmatrix}0&0\\1&0\end{smallmatrix} \bigr)-2z\phi(z)dz\otimes\Phi\bigl(\begin{smallmatrix}\frac{1}{2}&0\\0&-\frac{1}{2}\end{smallmatrix} \bigr)]\\
&=&d[z^2\phi(z)dz\otimes{\bigl(\begin{smallmatrix}0&-1&1\\1&0&0\\1&0&0\end{smallmatrix} \bigr)}-\phi(z)dz\otimes{\bigl(\begin{smallmatrix}0&1&1\\-1&0&0\\1&0&0\end{smallmatrix} \bigr)}-2z\phi(z)dz\otimes{\bigl(\begin{smallmatrix}0&0&0\\0&0&1\\0&1&0\end{smallmatrix} \bigr)}]\\
&& \text{by the definition of $\Phi$}\\
&=&d(z^2\phi(z)dz)\otimes{\bigl(\begin{smallmatrix}0&-1&1\\1&0&0\\1&0&0\end{smallmatrix} \bigr)}-d(\phi(z)dz)\otimes{\bigl(\begin{smallmatrix}0&1&1\\-1&0&0\\1&0&0\end{smallmatrix} \bigr)}-d(2z\phi(z)dz)\otimes{\bigl(\begin{smallmatrix}0&0&0\\0&0&1\\0&1&0\end{smallmatrix} \bigr)}\\
&=&0,~~~~~ \text{since $\phi(z)$ is holomorphic and hence $z^2\phi(z)$ and $2z\phi(z)$.}
\end{eqnarray*}
(b) We next prove that $\phi(z)dz\otimes\Phi\bigl(\begin{smallmatrix}-z&z^2\\-1&z\end{smallmatrix} \bigr)$ is coclosed, which is the heart of the lemma. 
From definition of $\delta$ (see equation (\ref{delta})), it is enough to show $d*(\sharp)(\phi(z)dz\otimes\Phi\bigl(\begin{smallmatrix}-z&z^2\\-1&z\end{smallmatrix} \bigr))=0$, which follows from
\begin{eqnarray}
&&d*(\sharp)[\phi(z)dz\otimes\Phi\bigl(\begin{smallmatrix}-z&z^2\\-1&z\end{smallmatrix} \bigr)]\nonumber\\
&=&d*\sharp[z^2\phi(z)dz\otimes\Phi\bigl(\begin{smallmatrix}0&1\\0&0\end{smallmatrix} \bigr)-\phi(z)dz\otimes\Phi\bigl(\begin{smallmatrix}0&0\\1&0\end{smallmatrix} \bigr)-2z\phi(z)dz\otimes\Phi\bigl(\begin{smallmatrix}\frac{1}{2}&0\\0&-\frac{1}{2}\end{smallmatrix} \bigr)]\nonumber\\
&=&d*[z^2\phi(z)dz\otimes\sharp(\Phi\bigl(\begin{smallmatrix}0&1\\0&0\end{smallmatrix} \bigr))-\phi(z)dz\otimes\sharp(\Phi\bigl(\begin{smallmatrix}0&0\\1&0\end{smallmatrix} \bigr))\nonumber\\
&&-2z\phi(z)dz\otimes\sharp(\Phi\bigl(\begin{smallmatrix}\frac{1}{2}&0\\0&-\frac{1}{2}\end{smallmatrix} \bigr))]\label{d*}
\end{eqnarray}
We then want to calculate $\sharp(\Phi\bigl(\begin{smallmatrix}0&1\\0&0\end{smallmatrix} \bigr)), \sharp(\Phi\bigl(\begin{smallmatrix}0&0\\1&0\end{smallmatrix} \bigr)),\sharp(\Phi\bigl(\begin{smallmatrix}\frac{1}{2}&0\\0&-\frac{1}{2}\end{smallmatrix} \bigr))$.

We choose a basis for $sl(3,\RR)$ as $\{E_1=\Phi\bigl(\begin{smallmatrix}0&1\\0&0\end{smallmatrix} \bigr)=\bigl(\begin{smallmatrix}0&-1&1\\1&0&0\\1&0&0\end{smallmatrix} \bigr),E_2=\Phi\bigl(\begin{smallmatrix}\frac{1}{2}&0\\0&-\frac{1}{2}\end{smallmatrix} \bigr)=\bigl(\begin{smallmatrix}0&0&0\\0&0&1\\0&1&0\end{smallmatrix} \bigr),E_3=\Phi\bigl(\begin{smallmatrix}0&0\\1&0\end{smallmatrix} \bigr)=\bigl(\begin{smallmatrix}0&1&1\\-1&0&0\\1&0&0\end{smallmatrix} \bigr),E_4=\bigl(\begin{smallmatrix}0&1&0\\0&0&0\\0&0&0\end{smallmatrix} \bigr),E_5=\bigl(\begin{smallmatrix}0&0&1\\0&0&0\\0&0&0\end{smallmatrix} \bigr),E_6=\bigl(\begin{smallmatrix}0&0&0\\0&1&0\\0&0&0\end{smallmatrix} \bigr),E_7=\bigl(\begin{smallmatrix}0&0&0\\0&1&0\\0&0&0\end{smallmatrix}\bigr),
E_8=\bigl(\begin{smallmatrix}0&0&0\\0&0&1\\0&0&0\end{smallmatrix}\bigr)\}$. The map $\sharp:sl(3,\RR)_{Ad\rho}\rightarrow sl(3,{\RR}^*)_{Ad{\rho}^*}$ is defined by setting
\begin{equation*}
\sharp(v)_x(u_x)=l_x(u_x,v_x),  ~~\text{for $u_x,v_x\in sl(3,{\RR}), x\in S$,}
\end{equation*}
we then have \begin{equation}
\sharp(A)_x=\sum_{\{1\leq i\leq 8\}}l(A,E_i)(E_i)^*,~~~ \text{for $A\in sl(3,\RR)$,}\label{sharp}
\end{equation}  where $(E_i)^*$ satisfies
$$(E_i)^*(E_j)=\left.\{ \begin{array}{rl}
1 &\mbox{ if $i=j$,} \\
0 &\mbox{ otherwise.}
       \end{array} \right. $$
Applying Lemma \ref{pairing g} to compute $l(E_i, E_j)$, for all $1\leq i\leq 3, 1\leq j\leq 8$, and then substitute the values into Equation (\ref{sharp}), we obtain the following formulas:
\begin{eqnarray*}
\sharp(E_1)=&\frac{4}{y^2}(E_1^*-x^2E_3^*+xE_2^*)+\frac{-1-x^2}{y^2}E_4^*+\frac{1-x^2}{y^2}E_5^*+\frac{2x}{y^2}E_8^*.\\
\sharp(E_2)=&\frac{4}{y^2}(xE_1^*-x(x^2+y^2)E_3^*+(x^2+\frac{1}{2}y^2)E_2^*)+\frac{-x^3-x-xy^2}{y^2}E_4^*\\&+\frac{-x^3+x-xy^2}{y^2}E_5^*+(\frac{2x^2}{y^2}+1)E_8^*.\\
\sharp(E_3)=&\frac{4}{y^2}(-x^2E_1^*+(x^2+y^2)^2E_3^*-x(x^2+y^2)E_2^*)
+\frac{x^4+y^4+2x^2y^2+x^2}{y^2}E_4^*\\
&+\frac{x^4+y^4+2x^2y^2-x^2}{y^2}E_5^*-\frac{2x}{y^2}E_8^*.
\end{eqnarray*}
Finally we apply the above formulas to compute Equation (\ref{d*})
\begin{eqnarray*}
&=&d*[z^2\phi(z)dz\otimes\sharp(E_1)-\phi(z)dz\otimes\sharp(E_3)-2z\phi(z)dz\otimes\sharp(E_2)]\\
&=&d*[-4\phi(z)dz\otimes E_1^*+4z^2\phi(z)dz\otimes E_3^*-4z\phi(z)dz\otimes E_2^*\\
&&+(z^2+1)\phi(z)dz\otimes E_4^*+(z^2-1)\phi(z)dz\otimes E_5^*-2z\phi(z)dz\otimes E_8^*]\\
&&\text{Observe that all terms inside the bracket remain holomorphic again.} \\
&=&d[-4i\overline{\phi(z)}d\overline{z}\otimes E_1^*+4i\overline{\phi(z)}{\overline{z}}^2d\overline{z}\otimes E_3^*-4i\overline{\phi(z)}\overline{z}d\overline{z}\otimes E_2^*\\
&&+i({\overline{z}}^2+1)\overline{\phi(z)}d\overline{z}\otimes E_4^*+i({\overline{z}}^2-1)\overline{\phi(z)}d\overline{z}\otimes E_5^*-2i\overline{z}\overline{\phi(z)}d\overline{z}\otimes E_8^*]\\
&=&0,~~~~~ \text{since $\phi(z)$ is holomorphic.}
\end{eqnarray*}
 
Thus (a) and (b) together imply that $\phi(z)dz\otimes\Phi\bigl(\begin{smallmatrix}-z&z^2\\-1&z\end{smallmatrix} \bigr)$ is harmonic. \qed
\end{subsection}\\

We now have the ingredients we need to prove Proposition \ref{Loftin metric}.\\

\textit{Proof of Proposition \ref{Loftin metric}}.  It is sufficient to prove the complexified version, i.e., we compute (with explanations of the steps given at the conclusion of the computation)
\begin{eqnarray*}
&&<[\phi(z)dz\otimes \Phi{\bigl(\begin{smallmatrix}-z&z^2\\-1&z\end{smallmatrix} \bigr)}],[\psi(z)dz\otimes \Phi{\bigl(\begin{smallmatrix}z&z^2\\-1&z\end{smallmatrix} \bigr)}]>_{Loftin}\\\\
&\myLongEqual{(i)}&Re\int_S (\phi(z)dz)\wedge *(\psi(z)dz)l_{f(z)}( \Phi{\bigl(\begin{smallmatrix}-z&z^2\\-1&z\end{smallmatrix} \bigr)},\Phi{\bigl(\begin{smallmatrix}z&z^2\\-1&z\end{smallmatrix} \bigr)})\\\\
&\myLongEqual{(ii)}&Re\int_S \phi(z)dz\wedge (i\overline{\psi(z)}d\overline{z})l_{f(z)}( \Phi{\bigl(\begin{smallmatrix}-z&z^2\\-1&z\end{smallmatrix} \bigr)},\Phi{\bigl(\begin{smallmatrix}z&z^2\\-1&z\end{smallmatrix} \bigr)})\\\\
&\myLongEqual{(iii)}&Re (16i\int_S\phi(z)\overline{\psi(z)}y^2 dzd\overline{z})\\\\
&\myLongEqual{(iv)}&32{<\phi dz^2,\psi dz^2>}_{WP}
\end{eqnarray*}
Therefore $\tau$ is an isometry of $\mathcal{B}(S)$ with the Loftin metric.\\
(i): Lemma \ref{harmonic} implies that $\phi(z)dz\otimes\Phi\bigl(\begin{smallmatrix}-z&z^2\\-1&z\end{smallmatrix} \bigr)$ is harmonic. Recall the definition of the Loftin metric, if $\sigma\otimes\phi$ is harmonic, we have \begin{equation*}
g_{\text{Loftin}}([\sigma\otimes \phi], [{\sigma}^{\prime}\otimes{\phi}^{\prime}])= \int_S(\sigma\wedge*{\sigma}{\prime})l(\phi, {\phi}^{\prime})dxdy.
\end{equation*}
(ii): Since $x, y$ are conformal coordinates for the Blaschke metric on the surface, we extend the action of the Hodge star operator to complex 1-forms by complex-antilinearity, i.e., $*(i\alpha)=-i*\overline{\alpha}$. From the definition of Hodge star (see equation (\ref{Hodge dual}), we see that $*dx=dy, *dy=-dx$, then $*{\phi(z)dz}=i\overline{\phi(z)}d\overline{z}$.\\
(iii): By Corollary \ref{pairing},  \begin{equation*}
l_{f(z)}( \Phi{\bigl(\begin{smallmatrix}-z&z^2\\-1&z\end{smallmatrix} \bigr)},\Phi{\bigl(\begin{smallmatrix}z&z^2\\-1&z\end{smallmatrix} \bigr)})=16y^{2}.
\end{equation*}
(iv): From the definition of the Weil-Petersson co-metric, ${<\phi dz^2,\psi dz^2>}_{WP}=Re\int_S \phi(z)\overline{\psi(z)}y^2 dxdy$ (see \cite{Goldman 1}, page 212) and $dz\wedge d\overline{z}=-2idx\wedge dy$. \qed
\end{subsection}
\end{subsection}
\end{section}
\begin{section}{Teichm$\ddot{\text{u}}$ller space is totally geodesic in $\mathcal{B}(S)$}\label{geodesic}
The goal of this section is to prove 
\begin{theorem}\label{total} Teichm$\ddot{\text{u}}$ller space endowed with (a constant multiple of) the Weil-Petersson metric is totally geodesic in $\mathcal{B}(S)$, endowed with the Loftin metric.
\end{theorem}

To achieve this goal, we make use of a dual map $\tau: \mathcal{B}(S)\rightarrow \mathcal{B}(S)$ which takes $(\Sigma, U)\rightarrow (\Sigma, -U)$, where $U$ is a holomorphic cubic differential on $\Sigma$. Therefore, the fixed set of this dual map $\tau$ is exactly the Teichm$\ddot{\text{u}}$ller locus $\cong \{(\Sigma,0)\}$. We will see it is sufficient to show the following theorem:
\begin{theorem}\label{Z} The dual map $\tau$ is an isometry of $\mathcal{B}(S)$ with respect to the Loftin metric.
\end{theorem} 
We now show how to derive Theorem \ref{total} from Theorem \ref{Z}.\\

\textit{\textit{Proof} ~of~ Theorem ~\ref{total}}. It is known (see \cite{Kabayashi}) that the fixed set of an isometry of a Riemannian manifold is a totally geodesic submanifold. Consider the manifold $\cB(S)$ endowed with the Loftin metric $g_{Loftin}$, we first have that the set $\{(\Sigma,0)\}\cong$ Teichm$\ddot{\text{u}}$ller space is the fixed set of the dual map $\tau$ on the manifold $\cB(S)$, next the dual map $\tau$ is an isometry of $\cB(S)$ from Theorem \ref{Z}, and moreover the Loftin metric restricts to be a constant multiple of the Weil-Petersson metric on the Teichm$\ddot{\text{u}}$ller locus from Theorem \ref{WP}. Combining with the fact stated in the first sentence,  we conclude that Teichm$\ddot{\text{u}}$ller space endowed with (a constant multiple of) the Weil-Petersson metric is totally geodesic in $\mathcal{B}(S)$ with the Loftin metric.\qed\\

So the remaining goal of this section is to show Theorem \ref{Z}. We divide the remaining part of this section into three parts. 

Part (I): we first define the conormal map of hyperbolic affine spheres and then state Proposition \ref{proposition} which tells us that the dual map $\tau$ is in fact induced by the conormal map $\nu$ of hyperbolic affine spheres. The reason we consider the conormal map $\nu$ instead of the dual map $\tau$ is that it is closely related to the definition of the Loftin metric, which involves the Cheng-Yau metric of the cone restricted to the hyperbolic affine sphere.

Part (II): we compare the Cheng-Yau metric on the cone restricted to the hyperbolic affine sphere $H$ and the one on the dual cone restricted to the dual hyperbolic affine sphere (image of the conormal map $\nu$) by showing Lemma \ref{h}, which will be essential in the proof of Theorem \ref{Z}.

Part (III): we describe the induced tangent map of the dual map $\tau$ on the tangent space of Goldman space $\cB(S)$. Then we continue to finish the proof of Theorem \ref{Z}.
\begin{subsection}*{Part (I)} We start with the definition of the conormal map $\nu$. Let $H\subset{\RR}^{3}$ be a nondegenerate hypersurface transverse to its position vector. Let ${\RR}_{3}$ be the dual space of ${\RR}^{3}$. We now define a map $\nu: H\rightarrow {\RR}_{3}$ as follows (see \S 5 in \cite{Nomizu}).

For each $p\in H$, let $\nu_p$ be the element of ${\RR}_{3}$ such that 
\begin{equation*}
\nu_p(\vec{p})=1
   ~~~~~\text{and}~~~~~ \nu_p(X)=0 ~~~~\text{for all}~~~~X\in T_p(H).
\end{equation*}
We have thus a differentiable map $\nu:H\rightarrow {\RR}_{3}$, called the conormal map. (This construction can be done with any transverse vector field $\xi$ in place of $\vec{p}$.)

\begin{proposition}(Schirokov-Schirokov \cite{SS}, unpublished work of Calabi, Gigena \cite{Gigena1} \cite{Gigena}) The image of the conormal map $\nu$ of a hyperbolic affine sphere $H$ with center 0 and affine mean curvature -1 is another such hyperbolic affine sphere $\overline{H}$ in the dual space ${\RR}_3$.
\end{proposition}
The conormal map can descend to be defined on the quotient of hyperbolic affine spheres (in other words, affine sphere structures), we obtain the following proposition which says that the dual map $\tau$ of Goldman space $\cB(S)$ is in fact induced by the conormal map $\nu$.
\begin{proposition}(Loftin \cite{Loftin 2}, \cite{Loftin 1})\label{proposition}  Given a properly convex $\RR P^2$-manifold $M=\Omega/\Gamma$, the conormal map $\nu$ with respect to the affine sphere structure induces a map to the dual manifold $M^*={\Omega}^*/{\Gamma}^*$, where ${\gamma}^*\in{\Gamma}^*$ is defined by ${\gamma}^*y(x)=y({\gamma}^{-1} x)$, for all $x\in {\RR}^3, y\in {\RR}_3, \gamma\in \Gamma$. This map is an isometry of the affine metrics. And the conormal map $\nu$ induces the dual map $\tau: \mathcal{B}(S)\rightarrow \mathcal{B}(S)$ takes $(\Sigma, U)\rightarrow (\Sigma, -U)$, where $U$ is a holomorphic cubic differential on $\Sigma$. 
\end{proposition}
\textit{Remark.} If we identify ${\RR}_3$ with ${\RR}^3$ by standard inner product, i.e., identify $y$ with $y^T$, we obtain an induced identification between $sl(3,\RR)$ and $sl(3,{\RR}^*)$, and between $SL(3,\RR)$ and $SL(3,{\RR}^*)$. Then we rephrase the description of ${\gamma}^*$ in the above proposition as follows:  for any $\gamma\in \Gamma$, we have
\begin{equation}
{\gamma}^*(y)=({\gamma}^T)^{-1}y, \text{for all $y\in {\RR}^3$}.\label{r}
\end{equation}

\begin{subsection}*{Part (II)} The goal of this part is to show the following lemma:
\begin{lemma}\label{h} Suppose $h$ is the matrix presentation of the Cheng-Yau metric (under the standard basis) on the cone restricted to the hyperbolic affine sphere inside the cone and $h^*$ is the one on the dual cone restricted to the dual hyperbolic affine sphere, then if we identify ${\RR}_3$ with ${\RR}^3$ by the standard inner product, we have $h_{\nu(p)}^*=h_{p}^{-1}$, for $p\in H$.
\end{lemma}
 \end{subsection}

\textit{Proof of Lemma \ref{h}}. Take a conformal parametrisation of the hyperbolic affine sphere $H$ as the map $f: \cD\rightarrow H$ with the coordinates $x,y$. Consider the conormal map $\nu: H\rightarrow \overline{H}\subset{\RR}_3$, composing with the map $f$, to obtain a map $\nu\circ f:\cD\rightarrow \overline {H}$. Abusing notation, we continue to write the new map $\nu\circ f$ as $\nu$, so that we then have 
\begin{equation}\label{33}
\nu(f)=1,\nu(f_x)=0,\nu(f_y)=0.
\end{equation} 
We next derive some properties of tangent map $\nu_*$:
\begin{lemma}\label{dual}(Proposition 5.1 in \cite{Nomizu}) Given a hyperbolic affine sphere $H$, the conormal map $\nu$ on $H$, and the affine metric $g=e^{\psi}|dz|^2$ along $H$, then
\begin{equation*}
\nu_{*}(Y)(\vec{p})=0 ~~~~   \text{and}   ~~~~ \nu_{*}(Y)(X)=-g(Y,X)~~~~
\text{for all}~~~~ X,Y\in T_p H.
\end{equation*}
\end{lemma} 
Applying Lemma \ref{dual}, and denote $\nu_x=\nu_{*}(\frac{\partial}{\partial x}), \nu_y=\nu_{*}(\frac{\partial}{\partial y}), f_x=f_{*}(\frac{\partial}{\partial x}), f_y=f_{*}(\frac{\partial}{\partial y})$, we have 
\begin{equation}\label{44}
\nu_x(f_x)=-g(f_x,f_x)=-e^{\psi},~~~\nu_x(f_y)=0,~~~ \nu_x(f)=0;
\end{equation}
\begin{equation}\label{55}
\nu_y(f_x)=0, ~~~\nu_y(f_y)=-g(f_y,f_y)=-e^{\psi},~~~\nu_y(f)=0.
\end{equation}\\
Collect equations (\ref{33}), (\ref{44}) and (\ref{55}) together, then we obtain\\
$\begin{pmatrix}
{\nu}\\
{-e^{-\frac{1}{2}\psi}\nu_x}\\
{-e^{-\frac{1}{2}\psi}\nu_y}
\end{pmatrix}\cdot (f, e^{-\frac{1}{2}\psi}f_x, e^{-\frac{1}{2}\psi}f_y)=
\begin{pmatrix}
1&0&0\\
0&1&0\\
0&0&1
\end{pmatrix}$.

Denote the matrix $A=(f,e^{-\frac{1}{2}\psi}f_x,e^{-\frac{1}{2}\psi}f_y)$, and we identify ${\RR}^3$ with ${\RR}_3$ by standard inner product, i.e., identify $v$ with $v^T$. Then, we have ${{A}^{-1}}^T=(\nu^T,-e^{-\frac{1}{2}\psi}\nu_x^T,{-e^{-\frac{1}{2}\psi}\nu_y^T})$.

We have the fact that pair $\{e^{-\frac{1}{2}\psi}f_x,e^{-\frac{1}{2}\psi}f_y\}$ is an orthonormal basis for the affine metric on $H$, and combining with Lemma \ref{2}, we then obtain that $\{f, e^{-\frac{1}{2}\psi}f_x,e^{-\frac{1}{2}\psi}f_y\}$ is an orthonormal basis for the Cheng-Yau metric on the cone restricted to $H$. Therefore we have \begin{equation}\label{11}
A^T h A=I.
\end{equation}

Similarily, $\{\nu^T,-e^{-\frac{1}{2}\psi}\nu_x^T,-e^{-\frac{1}{2}\psi}\nu_y^T\}$ is also an orthonormal basis for the Cheng-Yau metric on the cone ${\cC}^{*}\subset {\RR}_3$ restricted to the dual hyperbolic affine sphere $\overline{H}$, therefore we have
\begin{equation}\label{22}A^{-1} h^* {A^{-1}}^T=I.
\end{equation}
Combining equations (\ref{11}) and (\ref{22}), we obtain that $h^*=h^{-1}=AA^T$.\qed
\end{subsection}

\begin{subsection}*{Part (III)} We give an explicit description of  how the dual map $\tau$ affects the tangent vector by showing the following lemma.
\begin{lemma} \label{conormal}The dual map $\tau$ on $\mathcal{B}(S)$ induces the map ${\tau}_{*}$ on the tangent space $T_{[\rho]}\mathcal{B}(S)\cong H^1(S;sl(3,\RR)_{Ad \rho})$ as follows: \begin{eqnarray*} 
{\tau}_{*}: &H^1(S;sl(3,\RR)_{Ad\rho})&\rightarrow H^1(S;sl(3,{\RR})_{Ad{\rho}^*})\\
&[\sigma\otimes\phi]&\longmapsto [\sigma\otimes{\phi}^*]
\end{eqnarray*}
where $\forall s\in S$, we have 
\begin{equation}
\phi_s^*=-\phi_s^T, \text{where $[\sigma\otimes\phi]\in H^1(S;sl(3,\RR)_{Ad \rho})$}. \label{phi}
\end{equation}
\end{lemma}
\textit{Proof}. We carry out the proof in three steps.

Step 1: We first show that the tangent vector is represented by the cohomology class $[u]$ of 1-cocycle $u:\pi\rightarrow sl(3,\RR)$ satisfying $u(xy)-u(x)=Ad(\rho(x))(u(y))$, for $x, y\in \pi$. 

Consider a family $\{M_t\cong{\Omega}_t/{\Gamma}_t\}$ of convex $\RR P^2$-structures on the surface $S$, with corresponding family of conjugation classes of representations $\rho_t\in Hom(\pi, SL(3,\RR))$ (i.e., ${\Gamma}_t=\rho_t(\pi)$). 
Taking the derivative of both sides of the equation \begin{equation*}
\rho_t(xy)=\rho_t(x)\rho_t(y),~~~~ \text{for all $x, y\in \pi$},
\end{equation*} we obtain that the tangent vector at $\rho_0=\rho$ is a 1-cocycle $u:\pi\rightarrow sl(3,\RR)$ satisfying (see \cite{Goldman 1} for details) \begin{equation*}
u(xy)-u(x)=Ad(\rho(x))(u(y)).
\end{equation*}  

Step 2: We next want to show that the map  
\begin{equation} 
{\tau}_*: H^1(\pi,sl(3,{\RR}))\rightarrow H^1(\pi,sl(3,{\RR}))\nonumber\\
\end{equation} takes the tangent vector $u\in H^1(\pi,sl(3,{\RR}))$  to
\begin{equation} 
\tau_*(u)(\gamma)=-(u(\gamma))^T\in sl(3,\RR), \text{for all $\gamma\in \pi$}.\label{u}
\end{equation}
Considering the family $\tau(\rho_t)=\rho_t^*\in Hom(\pi, SL(3,\RR))$ and substituting $\gamma$ for ${\rho_t}^*(\gamma)$ into equation (\ref{r}), we obtain
\begin{equation}
{\rho_t}^*(\gamma)y=({\rho_t(\gamma)}^T)^{-1}y, \text{for all $y\in{\RR}^3$}\label{rho}
\end{equation} 
Taking the derivative of both sides of the above equation, we find that the tangent vector $\tau_*(u)$ satisfies
\begin{equation*}
\tau_*(u)(\gamma)y=-(u(\gamma))^Ty, \text{for all $y\in{\RR}^3$}
\end{equation*} 
Therefore, we obtain that \begin{equation}
\tau_*(u)(\gamma)=-(u(\gamma))^T\in sl(3,\RR),\text{ for all $\gamma\in \pi$.}\label{au}
\end{equation}

Step 3:  We want to show that $\tau_*([\sigma\otimes\phi])=[\sigma\otimes\phi^*]$.
Consider the canonical isomorphism from \begin{eqnarray}
&H^1(M,sl(3,{\RR}^*)_{Ad\rho})&\rightarrow H^1(\pi,sl(3,{\RR}))\nonumber\\
 &[\sigma\otimes\phi]&\longmapsto u_{\sigma\otimes\phi},\label{defn}
\end{eqnarray}
where  $u_{\sigma\otimes\phi}(\gamma):=\int_{\widetilde{\gamma}}\sigma\otimes\phi\in sl(3,\RR), \forall \gamma\in \pi$ and $\overline{\gamma}$ is an arbitrary closed curve in $S$ 
 representing $\gamma$. 

For any $\gamma\in \pi$, we compute  
\begin{eqnarray*}
u_{[\sigma\otimes\phi^*]}(\gamma)&=&\int_{\overline{\gamma}}\sigma\otimes\phi^*~~~ \text{by equation (\ref{defn})}\\
&=&-\int_{\overline{\gamma}}\sigma\otimes\phi^T~~~~\text{by equation (\ref{phi})}\\
&=&-(\int_{\overline{\gamma}}\sigma\otimes\phi)^T\\
&=&-(u_{[\sigma\otimes\phi]}(\gamma))^T ~~~ \text{by equation (\ref{defn})}\\
&=&\tau_*(u_{[\sigma\otimes\phi]})(\gamma)~~\text{after substituting $u_{[\sigma\otimes\phi]}$ for $u$ in equation (\ref{u})}.
\end{eqnarray*} 
Hence we have that $\tau_*(u_{[\sigma\otimes\phi]})=u_{[\sigma\otimes\phi^*]}$. Combining with  the injectivity of the map between $H^1(S,sl(3,\RR)_{Ad\rho})$ and $H^1(\pi,sl(3,\RR))$, we conclude that $\tau_*([\sigma\otimes\phi])=[\sigma\otimes\phi^*]$.
\qed\\

\textit{Proof of Theorem \ref{Z}}. We carry out the proof in three steps.

Step 1: We show that the Riemannian metrics on the two bundles \\$\mathcal{A}^1(S, sl(3,\RR)_{Ad{\rho}})$ and $\mathcal{A}^1(S, sl(3,\RR)_{Ad{\rho}^*})$ are isometric under $\mu: \sigma\otimes\phi\rightarrow\sigma\otimes{\phi}^*$.\\
Consider the Riemannian metric $l$ on the Lie algebra bundle $sl(3,\RR)_{Ad{\rho}}$, we compute
\begin{eqnarray*}
&&l(\phi,{\phi}^{\prime})|_x~~\text{taking $x\in S$}\\
&=&l_{f(\widetilde{x})}(A,B) ~~~~\text{by equation (\ref{l})}\\
&&\text{taking $\widetilde{\phi}_{f(\widetilde{x})}=A,\text{and} \widetilde{{\phi}^{\prime}}_{f(\widetilde{x})}=B$, where $\widetilde{x}$ is a preimage of $x$ in the domain $\cD\subset \CC$}\\
&=&tr(A^Th^{-1}Bh) ~~\text{by Lemma \ref{pairing g}}\\
&=&tr((A^Th^*)(B{h^{*}}^{-1}))~~ \text{by Lemma \ref{h}, $h^*=h^{-1}$} \\
&=&tr((B{h^{*}}^{-1})(A^Th^*))\\
&=&l^*_{\tau(\widetilde{x})}(B^T,A^T) ~~\text{by substituting $l$ with $l^*$ into Lemma \ref{pairing g}}\\
&=&l^*_{\tau(\widetilde{x})}(-A^T,-B^T)~~\text{by the symmetry and linearility of $l^*$}\\
&=&l^*({\phi}^*,{{\phi}^{\prime}}^*)|_x ~~ \text{by Lemma \ref{conormal}, $\widetilde{\phi}^*_{\tau(\widetilde{x})}=-A^T, \widetilde{{\phi}^{\prime}}^*_{\tau(\widetilde{x})}=-B^T$}
\end{eqnarray*} We denote this result as fact (1).

Furthermore, from Proposition \ref{proposition}, we see that the two Blaschke (or affine) metrics on $M\cong{\Omega}/{\Gamma}$ and $M^*\cong{\Omega}^*/{\Gamma}^*$ are isometric, hence we have fact (2):The Riemannian  metrics on space of 1-forms induced by the two Blaschke metrics are isometric. 
Hence we compute
\begin{eqnarray*}
&&g(\sigma\otimes {\phi}^*, {\sigma}^{\prime}\otimes{{\phi}^{\prime}}^*)\\
&=&\int_S(\sigma\wedge*{\sigma}^{\prime})l^*({\phi}^*, {{\phi}^{\prime}}^*)~~~\text{by definition of $g$ on $\mathcal{A}^1(S, sl(3,\RR)_{Ad{\rho}^*})$,}\\
&=&\int_S(\sigma\wedge*{\sigma}^{\prime})l(\phi, {\phi}^{\prime})~~~\text{by fact (1) and (2),}\\
&=&g(\sigma\otimes \phi, {\sigma}^{\prime}\otimes{\phi}^{\prime})~~~\text{by definition of $g$ on $\mathcal{A}^1(S, sl(3,\RR)_{Ad{\rho}})$.}
\end{eqnarray*}

Therefore the Riemannian metrics on the two bundles $\mathcal{A}^1(S, sl(3,\RR)_{Ad{\rho}})$ and $\mathcal{A}^1(S, sl(3,\RR)_{Ad{\rho}^*})$ are isometric under $\mu$. 

Step 2: We show that if $\sigma\otimes\phi$ is the harmonic representative in the cohomology class $ [\sigma\otimes\phi]$, then $\sigma\otimes{\phi}^*$ is also the unique harmonic representative in the cohomology class $ [\sigma\otimes{\phi}^*]$.

We begin by noting that, because $\sigma\otimes\phi$ is a harmonic representative, we have equivalently that, 
\begin{equation*}
d(\sigma\otimes\phi)=0, ~~~\delta(\sigma\otimes\phi)=0.
\end{equation*}

Note that $d$ is linear, $d(\sigma\otimes\phi)=0$ implies that $ d(-\sigma\otimes{\phi}^T)=0$, i.e., $ \sigma\otimes{\phi}^*=-\sigma\otimes{\phi}^T$ is closed.

Next we show that $ \sigma\otimes{\phi}^*=-\sigma\otimes{\phi}^T$ is coclosed. Now $\delta(\sigma\otimes\phi)=0$ and definition $\delta$ (see equation (\ref{delta})) implies that $d*\sharp(\sigma\otimes\phi)=0$, and thus $d*(\sigma\otimes\sharp{\phi})=0$. Suppose $\sharp$ and ${\sharp}^*$ are defined on the two bundles $\mathcal{A}^1(S, sl(3,\RR)_{Ad{\rho}})$ and $\mathcal{A}^1(S, sl(3,\RR)_{Ad{\rho}^*})$ respectively. Assume the basis of $sl(3,\RR)$ is $\{E_i^j\}$, i.e., where$\{E_i^j\}$ is $3\times3$ matrix whose entries are all 0 except for the $(i,j)$-entry which is equal to 1, where $1\leq i,j\leq 3$ and $i,j$ are not both 3. From the definition of $\sharp$ and ${\sharp}^*$, we have 
\begin{equation*}
\sharp(u)=l(u,E_i^j)(E_i^j)^*,~~~ {\sharp}^*(u)=l^*(u, E_i^j)(E_i^j)^*.
\end{equation*}
Hence $d*(\sigma\otimes\sharp{\phi})=0$ implies that $d*(\sigma\otimes l(\phi, E_i^j))(E_i^j)^*=0$, therefore 
\begin{equation}
d*(\sigma\otimes l(\phi, E_i^j))=0,~~\text{where $1\leq i,j\leq 3$ and $i,j$ are not both 3.} \label{d}
\end{equation}
Then we compute
\begin{eqnarray*}
&&d*(\sigma\otimes{\sharp}^*({-{\phi}^T}))\\
&=&d*(\sigma\otimes l^*(-{\phi}^T, E_i^j)(E_i^j)^*)\\
&=&d*(-\sigma\otimes l(\phi, E_j^i)(E_i^j)^*) \text{{   }Step 1 implies that  $l^*(-A^T,-B^T)=l(A,B)$}\\
&=&-d*(\sigma\otimes l(\phi, E_j^i))(E_i^j)^* \\
&=&0~~~~~{  }{ }{ } \text{by above equation (\ref{d}), as all the coefficients of $(E_i^j)^*$ vanish.} 
\end{eqnarray*}
Hence $d*{\sharp}^*(\sigma\otimes({-{\phi}^T}))=0$, and then $\delta(\sigma\otimes(-{\phi}^T))=0$. Therefore ${\tau}_*(\sigma\otimes{\phi})=\sigma\otimes(-{\phi}^T)$ is also a harmonic representative.

Step 3: Supposing $\sigma\otimes \phi$ is the harmonic representative in its cohomology class, we compute (explanations of each step are given in the end of the computations)
\begin{eqnarray*}
&&g_{\text{Loftin}}(\tau_*[\sigma\otimes\phi],\tau_*[{\sigma}^{\prime}\otimes{\phi}^{\prime}])\\
&\myLongEqual{(i)}&g_{\text{Loftin}}([\sigma\otimes {\phi}^*], {[\sigma}^{\prime}\otimes{{\phi}^{\prime}}^*])\\
&\myLongEqual{(ii)}&g(\sigma\otimes {\phi}^*, {\sigma}^{\prime}\otimes{{\phi}^{\prime}}^*)\\
&\myLongEqual{(iii)}&g(\sigma\otimes \phi, {\sigma}^{\prime}\otimes{\phi}^{\prime})\\
&\myLongEqual{(iv)}&g_{\text{Loftin}}([\sigma\otimes\phi],[{\sigma}^{\prime}\otimes{\phi}^{\prime}])
\end{eqnarray*}
(i): Lemma~~ \ref{conormal}~~ implies that $\tau_*( [\sigma\otimes\phi])= [\sigma\otimes\phi^*]$. \\
(ii): Because $\sigma\otimes \phi$ is the harmonic representative in its cohomology class, Step 2 tells us that $\sigma\otimes {\phi}^*$ is also the harmonic representative.\\
(iii): Step 1 implies that the Riemannian metrics on the two bundles $\mathcal{A}^1(S, sl(3,\RR)_{Ad{\rho}})$ and $\mathcal{A}^1(S, sl(3,\RR)_{Ad{\rho}^*})$ are isometric under $\mu$.\\
(iv): By definition of the Loftin metric by choosing harmonic represatatives.

Therefore, we conclude the theorem.
 \qed
\end{subsection}

\end{section}

\end{document}